\documentclass[3p]{elsarticle} 
\usepackage{enumitem}
\usepackage{placeins}
\usepackage{exscale}
\usepackage{mathrsfs}
\usepackage{relsize}
\usepackage{booktabs}
\usepackage{color}
\usepackage{multirow}
\usepackage{lineno,hyperref}
\usepackage{caption}
\usepackage{subfig}
\modulolinenumbers[5]

\usepackage{amsmath,amssymb,cases,amsthm}

\newtheorem{Definition}{Definition}[section]
\newtheorem{Theorem}{Theorem}[section]
\newtheorem{Lemma}{Lemma}[section]
\newtheorem{Proposition}{Proposition}[section]
\newtheorem{Remark}{Remark}[section]

\begin{document}

\begin{frontmatter}

\title{Conditional Mean-Field Type Linear Quadratic Optimal  Control Problem for Regime-Switching Diffusions with Indefinite Coefficients}
\tnotetext[mytitlenote]{This work is supported  by the National Key R\&D Program of China under Grant No. 2022YFA1006100,   the NSFC under Grant Nos. 61925306 and  12601823, and  the Shandong Postdoctoral Science Foundation.}

\author[]{Zhongbin Guo}
\ead{Guozhongbin0628@163.com}

\author[]{Guangchen Wang\corref{cor1}}
\ead{wguangchen@sdu.edu.cn}

\author[]{Shupeng Wang}
\ead{202590000013@sdu.edu.cn}

\address{School of Control Science and Engineering, Shandong University, Jinan 250061,   China}
\cortext[cor1]{Corresponding author.}

\begin{abstract}
We study a  conditional mean-field type linear quadratic optimal control problem for regime-switching diffusions with indefinite weights. The dynamics and cost depend on conditional expectations of the state and control given the history of the Markov chain. Under standard definiteness conditions, we establish the existence and uniqueness of the optimal control and derive its feedback representation and optimal cost through two sets of coupled  Riccati equations and a linear backward stochastic differential equation driven by the Markov chain. To treat the case of indefinite weights, we extend the equivalent cost functional approach to this setting using regime-dependent quadratic adjustments for the conditional mean of the state and its deviation from that mean. When the transformed weights in the equivalent cost functional satisfy standard definiteness conditions, this construction yields unique solvability of the original Hamiltonian system and a feedback representation of the unique optimal control. The proposed results also provide an effective method for establishing the solvability of fully coupled linear conditional mean-field forward-backward stochastic differential equations that do not satisfy certain existing monotonicity conditions. We further apply the results to conditional mean-variance portfolio selection with exogenous liability in an incomplete regime-switching market and derive the optimal investment strategy.
\end{abstract}


\end{frontmatter}

\section{Introduction}
A stochastic linear quadratic (LQ) optimal control problem refers to a problem in which  the cost functions have a quadratic form  in the state and control and the state is governed by a  linear stochastic differential equation (SDE). Owing to this elegant structure, one can derive an optimal control in state feedback form through a Riccati equation. Moreover, the solution of an LQ problem may provide a tractable  approximation to nonlinear control problems~\cite{AndersonMoore1989,Yong Zhou 1999}.  Therefore, it is of great importance to investigate the LQ problem. In a deterministic LQ problem, the control weight is usually assumed to be uniformly positive definite. A negative direction in this weight over a time interval makes the cost unbounded below. Early treatments of stochastic LQ control likewise imposed positive control weights~\cite{WonhamRiccati1968,BensoussanLectures1982}. Chen et al.~\cite{ChenLiZhou1998} showed that the stochastic LQ problem can nevertheless have a finite value and admit an optimal feedback control even when the control weight is indefinite. The distinction arises when the control enters the diffusion coefficient. In that case, the control affects the variance  of the state, and the resulting contribution can compensate for a negative control weight. 

A regime-switching diffusion (also referred to as a switching diffusion) couples a continuous solution of a SDE to a finite-state Markov chain that governs its drift and diffusion coefficients. A change of regime can represent a loss of production capacity after a machine failure or a shift in asset returns and volatilities as market conditions change~\cite{ElliottAggounMoore1995,Yin hybrid 2009}. For controlled switching diffusions, Donnelly~\cite{Donnelly2011} established a sufficient stochastic maximum principle. Subsequent work treated jump diffusions~\cite{ZhangElliottSiu2012}, forward-backward stochastic differential equations (FBSDEs) with delay~\cite{LvTaoWu2016}, and delayed jump diffusions~\cite{SavkuWeber2018}. In the LQ setting,  Zhang et al.~\cite{ZhangLiXiong2021} characterized open-loop optimality and closed-loop solvability  using  Hamiltonian system and the regular solution of the Riccati equations. Following \cite{ZhangLiXiong2021}, Wen et al.~\cite{WenLiXiongWeak} established the equivalence between open-loop and weak closed-loop solvability. For systems with random coefficients, Wen et al.~\cite{WenLiXiongZhang2023} proved solvability of the stochastic Riccati equation under uniform convexity and obtained a feedback representation of the optimal control.


This paper considers an LQ optimal control problem for regime-switching diffusions with two main features. First, the control weighting matrices are allowed to be indefinite, broadening the applicability of our results. Second, both the state equation and the cost functional depend on conditional expectations of the state and control given the history of the Markov chain. Here  we refer to such  problems  as conditional mean-field type control problems. Research on this class of problems began more recently, and fewer results are available. The inclusion of conditional expectation terms in our model is mainly inspired by Nguyen et al.~\cite{Nguyen George Yin Hoang 2020}. For a system of interacting particles that share the same Markov chain, Nguyen et al.~\cite{Nguyen George Yin Hoang 2020} showed that, as the number of particles tends to infinity, the empirical distribution of their states converges weakly to the law of a representative particle conditional on the history of the chain. Thus the limiting population distribution remains random even after individual fluctuations have been averaged out. This observation motivates conditional mean-field type control problems and games, in which the dynamics or costs depend on the conditional law of the state or control. In the LQ problem studied here, this dependence enters through the conditional expectations $\mathbb E[X_t|\mathcal F_t^\theta]$ and $\mathbb E[u_t|\mathcal F_t^\theta]$, where $\mathcal F_t^\theta$ represents the history of the Markov  chain up to time $t$.

For nonlinear conditional mean-field type control and game problems with regime switching, see \cite{Nguyen 2020,Nguyen general 2021,WuZhang2024,LiLvWuXiong2025,Guo Song Wang (2025)} for maximum principles and \cite{BayraktarCecchinChakraborty2023} for the dynamic programming principle.  For LQ models, Ding et al.~\cite{DingLiLvXu2025} derived and verified an optimal feedback control for an infinite-horizon discounted problem using two algebraic Riccati equations. Lv et al.~\cite{Lv Xiong Zhang 2023} obtained a feedback representation of a Stackelberg equilibrium, and Rol\'on Guti\'errez et al.~\cite{Rol贸n Guti茅rrez markovian 2024} used the solvability of conditional mean-field FBSDEs (CMF-FBSDEs) to establish the existence of a Nash equilibrium. Related results include  a separation principle for LQ control under partial observation~\cite{Guo Wang 2026}. These control and game formulations~\cite{DingLiLvXu2025,Lv Xiong Zhang 2023,Rol贸n Guti茅rrez markovian 2024,Guo Wang 2026} assume positive definite control weights. 

This paper studies a conditional mean-field type LQ problem for regime-switching diffusions with indefinite coefficients. To address the indefinite coefficients, we use the equivalent cost functional method introduced by Yu~\cite{Yu2013} and further developed by Huang and Yu~\cite{HuangYu2014}. The core idea of this method is to construct a new cost functional that is equivalent to the original one and satisfies the standard definiteness conditions. Here, we briefly discuss the differences between our work and the existing studies~\cite{LiLiYu2020,WangWang2024,MeiWeiYong2024}. Li et al.~\cite{LiLiYu2020} studied an indefinite mean-field type LQ problem and introduced the concept of a relaxed compensator by developing the equivalent cost functional method. Wang and Wang~\cite{WangWang2024} subsequently extended the results in~\cite{LiLiYu2020} to jump-diffusion models. The studies~\cite{LiLiYu2020,WangWang2024} use unconditional expectations without Markovian regime switching, whereas our work considers conditional mean-field type LQ problems with regime switching. Here, the conditional expectation of the state is a random process adapted to the filtration generated by the Markov chain. The switching process enables our model to describe more complex phenomena but also makes the analysis more involved. Since the coefficients are modulated by the switching process, we use two regime-dependent matrix families, \(H\) and \(K\), to construct the equivalent cost functional. They are associated with the deviation of the state from its conditional expectation and the conditional expectation itself, respectively. The regime-switching mechanism introduces coupling across regimes in both the transformed cost weights and the Riccati equations. We establish the corresponding transformations of the Hamiltonian systems and coupled Riccati equations, together with the associated linear BSDE driven by the Markov chain. Corrections involving differences between the adjustment matrices across regimes also arise in the transformation of the jump integrands in the adjoint equation. Mei et al.~\cite{MeiWeiYong2024} investigated an indefinite conditional mean-field type LQ problem that includes regime switching and characterized finite-horizon open-loop and closed-loop solvability. Although our problem setting is similar to that of~\cite{MeiWeiYong2024}, our approach and solvability results differ in essential respects. Mei et al.~\cite{MeiWeiYong2024} proved that uniform convexity of the cost functional is equivalent to strongly regular solvability of the Riccati equations. Our approach complements this characterization with a sufficient criterion expressed as explicit inequalities for the transformed cost coefficients. When suitable matrices \(H\) and \(K\) make these coefficients satisfy the standard definiteness conditions, we obtain unique solvability of the original Hamiltonian system and coupled Riccati equations. These equations yield the unique optimal control in feedback form and its optimal cost. The auxiliary matrices need only satisfy the coefficient inequalities and need not solve the original Riccati equations. In Example~2, these inequalities are verified directly using constant choices of \(H\) and \(K\). This establishes unique solvability of a fully coupled CMF-FBSDE that fails the specific monotonicity conditions examined there.

The main contributions of this paper are summarized as follows.
\begin{itemize}
	\item As a foundation, we first consider the case in which the standard definiteness conditions hold. We prove the existence and uniqueness of the optimal control. The  feedback representation and the optimal cost are derived through two sets of coupled  Riccati equations and a linear BSDE driven by the Markov chain. The unique solvability of the involved Riccati equations is also established.

	\item For indefinite coefficients, we establish a sufficient solvability criterion through explicit inequalities for the transformed cost weights. The criterion yields unique solvability of the Hamiltonian system and coupled Riccati equations, as well as the optimal feedback control and its cost. Conversely, Riccati solutions with uniformly positive definite effective control weights provide the required adjustment matrices.
	
	\item We provide two examples and a financial application to illustrate the effectiveness of our results. In Example 1, we derive an explicit optimal control and optimal cost for an LQ problem with negative control weights using the equivalent cost functional method. In Example 2, we demonstrate how this method can establish the global solvability of fully coupled CMF-FBSDEs through their connection with the Hamiltonian systems of indefinite LQ control problems. We apply this approach to prove the solvability of a system that does not satisfy the solvability conditions considered in the existing literature. In the financial application, we solve a conditional mean-variance portfolio selection problem with exogenous liability in an incomplete regime-switching market, deriving the optimal investment strategy and the relation between the conditional mean and variance of optimal terminal net wealth.
\end{itemize}

The remainder of this paper is organized as follows. Section~\ref{sec: Problem Formulation} gives the problem formulation. Section~\ref{sec: posi def case} studies the problem under standard definiteness conditions. Section~\ref{sec: CMF of indef coe}  investigates   the indefinite  coefficients case by using the equivalent cost functional approach. Section~\ref{sec:example} presents two illustrative examples. Section~\ref{sec: fo alm} applies the results to conditional mean-variance portfolio selection with liability and provides numerical illustrations. Finally, Section \ref{sec: con rema} concludes the paper and outlines directions for future research.

\section{Problem Formulation}
\label{sec: Problem Formulation}
For positive integers $n$ and $m$, let $\mathbb R^n$, $\mathbb R^{n\times m}$, and $\mathbb S^n$ denote the spaces of real $n$-vectors, $n\times m$ matrices, and symmetric $n\times n$ matrices, respectively. We write  $M^\top$ for the transpose of the matrix $M$. For vectors, $\langle\cdot,\cdot\rangle$ denotes the Euclidean inner product and $|\cdot|$ denotes the corresponding norm. For $x\in\mathbb R^n$ and $M\in\mathbb S^n$, let $|x|_M^2=x^\top Mx$ denote the quadratic form. For $Q,R\in\mathbb S^n$, write $Q\geq0$ if $Q$ is positive semidefinite and $Q\geq R$ if $Q-R\geq0$. Let $\mathbb I_{\mathcal A}$ denote the indicator of a measurable set $\mathcal A$. The identity matrix is denoted by $\mathbf I_{n\times n}$. $\mathbf 1_m$ denotes the $m$-vector of ones.

Fix a finite horizon $T>0$, and let $\mathbb N=\{1,\ldots,N\}$ and $\mathbb D=\{1,\ldots,d\}$. On a probability space $(\Omega,\mathcal F,\mathbb P)$, let $W=(W^1,\ldots,W^N)$ be an $N$-dimensional standard Brownian motion, $\theta$ a time-homogeneous Markov chain with state space $\mathbb D$, and $\xi$ an $\mathbb R^n$-valued square-integrable random variable. Assume that $W$, $\theta$, and $\xi$ are mutually independent. Write $\ell_0:=\theta_0$. Let $\mathbb F$ and $\mathbb F^\theta$ be the usual augmentations of the filtrations generated by $(\xi,W,\theta)$ and $\theta$, respectively. For a filtration $\mathbb F^\beta$ used below, write $\mathbb E_t^\beta[\cdot]=\mathbb E[\cdot\mid\mathcal F_t^\beta]$. Let $\Lambda=(\lambda^{\ell j})_{\ell,j\in\mathbb D}$ be the generator of $\theta$, where $\lambda^{\ell j}\geq0$ for $\ell\ne j$ and $\lambda^{\ell\ell}=-\sum_{j\ne\ell}\lambda^{\ell j}$. For $\ell\ne j$, define the transition count and its compensated process by
\(
N_t^{\ell j}=\sum_{0<s\leq t}\mathbb I_{\{\theta_{s-}=\ell,\,\theta_s=j\}}\) and \(
M_t^{\ell j}=N_t^{\ell j}-\lambda^{\ell j}\int_0^t\mathbb I_{\{\theta_{s-}=\ell\}}\,\mathrm ds.
\)
Then $M^{\ell j}$ is a purely discontinuous, square-integrable $\mathbb F^\theta$-martingale \cite[Lemma IV.21.12]{Rogers 2000}, with optional and predictable quadratic variations given by 
\(
[M^{\ell j}]_t=N_t^{\ell j}\) and 
\(\langle M^{\ell j}\rangle_t=\lambda^{\ell j}\int_0^t\mathbb I_{\{\theta_{s-}=\ell\}}\,\mathrm ds.
\)
Set $M^{\ell\ell}=[M^{\ell\ell}]=\langle M^{\ell\ell}\rangle=0$. Moreover, for $i\in\mathbb N$ and distinct pairs $(\ell,j),(k,\iota)\in\mathbb D\times\mathbb D$, the following orthogonality relations hold \cite[Section 1.8]{Lipster 1989}:
\(
[W^i,M^{\ell j}]=0\) and \([M^{\ell j},M^{k\iota}]=0.
\)

For a finite-dimensional Euclidean space $\mathbb H$ and a filtration $\mathbb G$ satisfying the usual conditions, we use the following spaces. $\mathcal L^\infty([0,T];\mathbb H)$ is the space of Borel-measurable functions $\mathcal X\colon[0,T]\to\mathbb H$ that are uniformly bounded on $[0,T]$.  $\mathcal S_{\mathbb G}^2([0,T];\mathbb H)$ is the space of $\mathbb G$-adapted c\`adl\`ag processes $\mathcal X\colon[0,T]\times\Omega\to\mathbb H$ such that $\mathbb E[\sup_{t\in[0,T]}|\mathcal X_t|^2]<\infty$. $\mathcal L_{\mathbb G}^2([0,T];\mathbb H)$ is the space of $\mathbb G$-progressively measurable processes $\mathcal X\colon[0,T]\times\Omega\to\mathbb H$ such that $\mathbb E[\int_0^T|\mathcal X_t|^2\,\mathrm dt]<\infty$.  $\mathcal M_{\mathbb G}^2([0,T];\mathbb H)$ is the space of families $\mathcal X=\{\mathcal X^{\ell j}\}_{\ell,j\in\mathbb D}$ of $\mathbb G$-predictable processes $\mathcal X^{\ell j}\colon[0,T]\times\Omega\to\mathbb H$ such that $\mathcal X^{\ell\ell}=0$ for every $\ell\in\mathbb D$ and
\(
\sum_{\ell,j=1}^{d}\mathbb E\![\int_0^T|\mathcal X_t^{\ell j}|^2\,\mathrm d\langle M^{\ell j}\rangle_t]<\infty.
\)  For $\mathcal X\in\mathcal M_{\mathbb G}^2([0,T];\mathbb H)$, we use the notation
\(
\int_0^t\mathcal X_s\bullet\mathrm dM_s
:=\sum_{\ell,j=1}^{d}\int_0^t\mathcal X_s^{\ell j}\,\mathrm dM_s^{\ell j}\) and \(
\mathcal X_t\bullet\mathrm dM_t
:=\sum_{\ell,j=1}^{d}\mathcal X_t^{\ell j}\,\mathrm dM_t^{\ell j}.
\)

Consider the conditional mean-field SDE  with regime switching
\begin{equation}\label{equ: com sta}
	\left\{\begin{aligned}
		\mathrm{d}X^{u}_t=&\ \bigl[A_{(\theta_t,t)}X^{u}_t+\widehat{A}_{(\theta_t,t)}\mathbb{E}^{\theta}_t\bigl[X^{u}_t\bigr]+B_{(\theta_t,t)}u_{t}+\widehat{B}_{(\theta_t,t)}\mathbb{E}^{\theta}_t\bigl[u_{t}\bigr]+\alpha_{t}\bigr]\mathrm{d}t\\
		&\  +\sum_{i=1}^{N}\Bigl[C^{i}_{(\theta_t,t)}X^{u}_t+\widehat{C}^{i}_{(\theta_t,t)}\mathbb{E}^{\theta}_t\bigl[X^{u}_t\bigr]+D^{i}_{(\theta_t,t)}u_{t}+\widehat{D}^{i}_{(\theta_t,t)}\mathbb{E}^{\theta}_t\bigl[u_{t}\bigr]+\beta^{i}_{t}\Bigr]\mathrm{d}W^i_t,\\
		X^{u}_0=&\ \xi,
	\end{aligned}\right.
\end{equation}
where $X^u$ is the $\mathbb R^n$-valued state process and $u$ is the $\mathbb R^m$-valued control process. We impose the following assumption on the coefficients.
\begin{enumerate}[label=$\mathbf{(A\arabic*)}$,leftmargin=3em]
	\item\label{ass:coef in com dyna} For every $\ell\in\mathbb D$ and $i\in\mathbb N$, \(A_{(\ell,\cdot)},\widehat A_{(\ell,\cdot)},C^i_{(\ell,\cdot)},\widehat C^i_{(\ell,\cdot)}\in\mathcal L^\infty([0,T];\mathbb R^{n\times n})\) and \(B_{(\ell,\cdot)},\widehat B_{(\ell,\cdot)},D^i_{(\ell,\cdot)},\widehat D^i_{(\ell,\cdot)}\in\mathcal L^\infty([0,T];\mathbb R^{n\times m})\). Moreover, $\alpha,\beta^i\in\mathcal L_{\mathbb F^\theta}^2([0,T];\mathbb R^n)$ for every $i\in\mathbb N$.
\end{enumerate}

The admissible control set is
\(
\mathcal U:=\mathcal L_{\mathbb F}^2([0,T];\mathbb R^m).
\)
Under assumption \ref{ass:coef in com dyna}, for every $u\in\mathcal U$, the state equation \eqref{equ: com sta} admits a unique solution $X^u\in\mathcal S_{\mathbb F}^2([0,T];\mathbb R^n)$; see Lemma 2.1 of \cite{Nguyen 2020}. The cost functional is
\begin{equation}\label{equ: com cos}
	\begin{aligned}
		J\bigl(\xi,\ell_0; u\bigr)
		= &\  \frac{1}{2} \mathbb{E}\biggl\{\int _ { 0 } ^ { T } \biggl[\left\langle\left[\begin{array}{l}
			X^{u}_t \\
			u_{t}
		\end{array}\right],\left[\begin{array}{cc}
			Q_{(\theta_t,t)} & S_{(\theta_t,t)} \\
			S^{\top}_{(\theta_t,t)} & R_{(\theta_t,t)}
		\end{array}\right]\left[\begin{array}{l}
			X^{u}_t \\
			u_{t}
		\end{array}\right]\right\rangle+2\left\langle\left[\begin{array}{l}
			X^{u}_t \\
			u_{t}
		\end{array}\right],\left[\begin{array}{l}
			q_{t} \\
			r_{t}
		\end{array}\right]\right\rangle\\
		&\ +\left\langle\left[\begin{array}{l}
			\mathbb{E}^{\theta}_t\bigl[X^{u}_t\bigr] \\
			\mathbb{E}^{\theta}_t\bigl[u_{t}\bigr]
		\end{array}\right],\left[\begin{array}{cc}
			\widehat{Q}_{(\theta_t,t)} & \widehat{S}_{(\theta_t,t)} \\
			\widehat{S}^{\top}_{(\theta_t,t)} & \widehat{R}_{(\theta_t,t)}
		\end{array}\right]\left[\begin{array}{c}
			\mathbb{E}^{\theta}_t\bigl[X^{u}_t\bigr] \\
			\mathbb{E}^{\theta}_t\bigl[u_{t}\bigr]
		\end{array}\right]\right\rangle +2\left\langle\left[\begin{array}{l}
			\mathbb{E}^{\theta}_t\bigl[X^{u}_t\bigr] \\
			\mathbb{E}^{\theta}_t\bigl[u_{t}\bigr]
		\end{array}\right],\left[\begin{array}{l}
			\widehat{q}_{t} \\
			\widehat{r}_{t}
		\end{array}\right]\right\rangle\biggr]\mathrm{d} t\\
		&\  +\langle X^{u}_T, G_{\theta_T} X^{u}_T\rangle+2\langle g, X^{u}_T\rangle+\langle\mathbb{E}^{\theta}_T\bigl[X^{u}_T\bigr], \widehat{G}_{\theta_T} \mathbb{E}^{\theta}_T\bigl[X^{u}_T\bigr]\rangle+2\langle \widehat{g}, \mathbb{E}^{\theta}_T\bigl[X^{u}_T\bigr]\rangle\biggr\}.
	\end{aligned}
\end{equation}
We impose the following assumption on its coefficients.
\begin{enumerate}[label=$\mathbf{(A\arabic*)}$,leftmargin=3em,start=2]
	\item\label{ass: coe in com cost} For every $\ell\in\mathbb D$, \(Q_{(\ell,\cdot)},\widehat Q_{(\ell,\cdot)}\in\mathcal L^\infty([0,T];\mathbb S^n)\), \(S_{(\ell,\cdot)},\widehat S_{(\ell,\cdot)}\in\mathcal L^\infty([0,T];\mathbb R^{n\times m})\), \(R_{(\ell,\cdot)},\widehat R_{(\ell,\cdot)}\in\mathcal L^\infty([0,T];\mathbb S^m)\), and \(G_\ell,\widehat G_\ell\in\mathbb S^n\). The matrices $G_\ell$ and $\widehat G_\ell$ are deterministic. Furthermore,
	\(
	q,\widehat q\in\mathcal L_{\mathbb F^\theta}^2([0,T];\mathbb R^n)\), 
	\(r,\widehat r\in\mathcal L_{\mathbb F^\theta}^2([0,T];\mathbb R^m)
	\), 
	and $g,\widehat g$ are square-integrable, $\mathcal F_T^\theta$-measurable $\mathbb R^n$-valued random variables.
\end{enumerate}


The conditional mean-field type LQ problem considered in this paper is formulated as follows.\\
\noindent\textbf{Problem CMF.} Find an admissible control $u^*\in\mathcal U$ such that
\[
J\bigl(\xi,\ell_0;u^*\bigr)=\inf_{u\in\mathcal U}J\bigl(\xi,\ell_0;u\bigr),
\]
subject to the state equation \eqref{equ: com sta}.

We first study \textbf{Problem CMF} under the positive definiteness assumption in Section \ref{sec: posi def case}. In Section \ref{sec: CMF of indef coe}, we employ the equivalent cost functional method to solve \textbf{Problem CMF} in the case of indefinite coefficients.

Throughout the paper, set $\widetilde\eta=\eta+\widehat\eta$, where $\eta\in\{A,B,C^i,D^i,S,Q,R,G,q,r,g\}$. When no confusion can arise, we omit the argument $(\theta_t,t)$ from regime-dependent coefficients evaluated at $(\theta_t,t)$.

\section{Positive Definite Coefficients Case}\label{sec: posi def case}
First, we introduce the  standard positive definiteness assumption.
\begin{enumerate}[label=$\mathbf{(PD)}$,leftmargin=3em]
	\item\label{ass: stand cod} There exists a constant $\delta>0$ such that, for every $\ell\in\mathbb{D}$ and $t\in[0,T]$, $R_{(\ell,t)},\widetilde{R}_{(\ell,t)}\geq \delta \mathbf{I}_{m\times m}$, $Q_{(\ell,t)}-S_{(\ell,t)}R^{-1}_{(\ell,t)}S^{\top}_{(\ell,t)},\widetilde{Q}_{(\ell,t)}-\widetilde{S}_{(\ell,t)}\widetilde{R}^{-1}_{(\ell,t)}\widetilde{S}^{\top}_{(\ell,t)}\geq 0$, and $G_{\ell},\widetilde{G}_{\ell}\geq 0$.
\end{enumerate}

The following result establishes the existence and uniqueness of the optimal control for \textbf{Problem CMF}.
\begin{Lemma}\label{lem: unique mini}
	Under assumptions \ref{ass:coef in com dyna}, \ref{ass: coe in com cost}, and \ref{ass: stand cod}, \textbf{Problem CMF} admits a unique optimal control.
\end{Lemma}
\begin{proof}
	The proof is provided in  \ref{app: proof unique mini}.
\end{proof}

We next give a necessary and sufficient condition for optimality in \textbf{Problem CMF}.
\begin{Lemma}\label{lem: nece and suff of optim}
	Under assumptions \ref{ass:coef in com dyna}, \ref{ass: coe in com cost}, and \ref{ass: stand cod}, $(X^{*},u^{*})$ is an optimal pair for \textbf{Problem CMF} if and only if there exist processes $(p,\zeta,k)$ such that $(X^{*},u^{*},p,\zeta,k)$ satisfies the following Hamiltonian system:
	\begin{equation}\label{equ: Ham for com}
		\left\{\begin{aligned}
			\mathrm{d}X^{*}_{t}=&\ \bigl[AX^{*}_{t}+\widehat{A}\mathbb{E}^{\theta}_t\bigl[X^{*}_{t}\bigr]+Bu^{*}_{t}+\widehat{B}\mathbb{E}^{\theta}_t\bigl[u^{*}_{t}\bigr]+\alpha_{t}\bigr]\mathrm{d}t\\ &\ +\sum_{i=1}^{N}\Bigl[C^{i}X^{*}_{t}+\widehat{C}^{i}\mathbb{E}^{\theta}_t\bigl[X^{*}_{t}\bigr]+D^{i}u^{*}_{t}+\widehat{D}^{i}\mathbb{E}^{\theta}_t\bigl[u^{*}_{t}\bigr]+\beta^{i}_{t}\Bigr]\mathrm{d}W^i_t,\\
			\mathrm{d}p_t=&\   -\Bigl[A^{\top}p_t+\widehat{A}^{\top}\mathbb{E}^{\theta}_t\bigl[p_t\bigr]+\sum_{i=1}^{N}C^{i\top}\zeta^{i}_t+\sum_{i=1}^{N}\widehat{C}^{i\top}\mathbb{E}^{\theta}_t\bigl[\zeta^{i}_t\bigr]+QX^{*}_{t} +\widehat{Q}\mathbb{E}^{\theta}_t\bigl[X^{*}_{t}\bigr]\\
			&\  +Su^{*}_{t} +\widehat{S}\mathbb{E}^{\theta}_t\bigl[u^{*}_{t}\bigr]+\widetilde{q}_{t}\Bigr]\mathrm{d}t +\sum_{i=1}^{N}\zeta^{i}_t\mathrm{d}W^{i}_t+k_t\bullet\mathrm{d}M_t,\\
			\quad X^{*}_0=&\ \  \xi,\quad  p_T=\ G_{\theta_T}X^{*}_T+\widehat{G}_{\theta_T}\mathbb{E}^{\theta}_T\bigl[X^{*}_T\bigr]+\widetilde{g},\\
			&\ \hspace{-1cm} B^{\top}p_t+\widehat{B}^{\top}\mathbb{E}^{\theta}_t\bigl[p_t\bigr]+\sum_{i=1}^{N}D^{i\top}\zeta^{i}_t+\sum_{i=1}^{N}\widehat{D}^{i\top}\mathbb{E}^{\theta}_t\bigl[\zeta^{i}_t\bigr] +S^{\top}X^{*}_{t}+\widehat{S}^{\top}\mathbb{E}^{\theta}_t\bigl[X^{*}_{t}\bigr]+Ru^{*}_{t}+\widehat{R}\mathbb{E}^{\theta}_t\bigl[u^{*}_{t}\bigr]+\widetilde{r}_{t}=0.
		\end{aligned}\right.
	\end{equation}
\end{Lemma}
\begin{proof}
	This result gives a maximum principle for \textbf{Problem CMF}, with the last equality in \eqref{equ: Ham for com} expressing the maximum condition. The proof uses standard convex variational techniques and a duality argument, following the approach of Theorem~3.7 in \cite{Nguyen 2020}. We omit the details.
\end{proof}

\begin{Lemma}\label{lem: well-pose of FBSDE}
	Under assumptions \ref{ass:coef in com dyna}, \ref{ass: coe in com cost}, and \ref{ass: stand cod}, the FBSDE \eqref{equ: Ham for com} admits a unique solution $(X^{*},u^{*},p,$ $\zeta,k)\in\mathcal S_{\mathbb F}^2([0,T];\mathbb R^n)\times \mathcal L_{\mathbb F}^2([0,T];\mathbb R^m)\times\mathcal S_{\mathbb F}^2([0,T];\mathbb R^n)\times (\mathcal L_{\mathbb F}^2([0,T];\mathbb R^n))^{N}\times \mathcal{M}^{2}_{\mathbb{F}}([0,T];\mathbb{R}^{n})$.
\end{Lemma}
\begin{proof}
	Lemma~\ref{lem: unique mini} and the well-posedness of \eqref{equ: com sta} give a unique optimal pair $(X^*,u^*)$. By Lemma~\ref{lem: nece and suff of optim}, there exist processes $(p,\zeta,k)$ such that $(X^*,u^*,p,\zeta,k)$ solves \eqref{equ: Ham for com}. Let $(\bar X,\bar u,\bar p,\bar\zeta,\bar k)$ be another solution of \eqref{equ: Ham for com}. Lemma~\ref{lem: nece and suff of optim} implies that $(\bar X,\bar u)$ is also optimal. Hence, $\bar u=u^*$ in $\mathcal U$ and $\bar X=X^*$. For the fixed pair $(X^*,u^*)$, the uniqueness of solutions to conditional mean-field BSDEs with regime switching (see Theorem 3.4 in \cite{Nguyen 2020}) yields that $p=\bar p$, $\zeta=\bar\zeta$, and $k^{\ell j}=\bar k^{\ell j}$  for every $\ell,j\in\mathbb D$. This proves uniqueness.
\end{proof}

We now decouple the Hamiltonian system \eqref{equ: Ham for com} using the four-step scheme  approach and construct the optimal
control in feedback form. To this end, we introduce  two  families of coupled Riccati equations
\begin{equation}\label{equ: Ricca one of Pro CMF}
	\left\{\begin{aligned}
		\dot{\Phi}_{(\ell,t)}=&\ -\Bigl[A^{\top}_{(\ell,t)}\Phi_{(\ell,t)}+\Phi_{(\ell,t)} A_{(\ell,t)}+\sum_{i=1}^{N}C^{i\top}_{(\ell,t)}\Phi_{(\ell,t)} C^{i}_{(\ell,t)}+Q_{(\ell,t)}-L_{(\ell,t)}\Sigma^{-1}_{(\ell,t)}L^{\top}_{(\ell,t)}+\sum_{j\in\mathbb{D}}\lambda^{\ell j}\Phi_{(j,t)}\Bigr],\\
		\Phi_{(\ell,T)}=&\ G_{\ell},
	\end{aligned}\right.
\end{equation}
\begin{equation}\label{equ: Ricca two of Pro CMF}
	\left\{\begin{aligned}
		\dot{\Pi}_{(\ell,t)}=&\ -\Bigl[\widetilde{A}^{\top}_{(\ell,t)}\Pi_{(\ell,t)}+\Pi_{(\ell,t)} \widetilde{A}_{(\ell,t)}+\sum_{i=1}^{N}\widetilde{C}^{i\top}_{(\ell,t)}\Phi_{(\ell,t)} \widetilde{C}^{i}_{(\ell,t)}+\widetilde{Q}_{(\ell,t)}-\widetilde{L}_{(\ell,t)}\widetilde{\Sigma}^{-1}_{(\ell,t)}\widetilde{L}^{\top}_{(\ell,t)}+\sum_{j\in\mathbb{D}}\lambda^{\ell j}\Pi_{(j,t)}\Bigr],\\
		\Pi_{(\ell,T)}=&\ \widetilde{G}_{\ell},
	\end{aligned}\right.
\end{equation}and a linear regime-switching BSDE
\begin{equation}\label{equ: Lya of Pro CMF}
	\left\{\begin{aligned}
		\mathrm{d}\varphi_{t}=&\ -\Bigl[\bigl(\widetilde{A}-\widetilde{B}\widetilde{\Sigma}^{-1}\widetilde{L}^{\top}\bigr)^{\top}\varphi_{t}+\sum_{i=1}^{N}\bigl(\widetilde{C}^{i}-\widetilde{D}^{i}\widetilde{\Sigma}^{-1}\widetilde{L}^{\top}\bigr)^{\top}\Phi\beta^{i}_{t}-\widetilde{L}\widetilde{\Sigma}^{-1} \widetilde{r}_{t}+\widetilde{q}_{t}+\Pi\alpha_{t}\Bigr]\mathrm{d}t+\rho_t\bullet\mathrm{d}M_t,\\
		\varphi_{T}=&\ \widetilde{g},
	\end{aligned}\right.
\end{equation}where 
\begin{equation}\label{equ: nota Sigma}
	\left\{\begin{aligned}
		\Sigma_{(\ell,t)}=&\  R_{(\ell,t)}+\sum_{i=1}^{N}D^{i\top}_{(\ell,t)}\Phi_{(\ell,t)}D^{i}_{(\ell,t)}, \quad&\  L_{(\ell,t)}=&\   S_{(\ell,t)}+\Phi_{(\ell,t)} B_{(\ell,t)}+\sum_{i=1}^{N}C^{i\top}_{(\ell,t)}\Phi_{(\ell,t)} D^{i}_{(\ell,t)}, \\
		\widetilde{\Sigma}_{(\ell,t)}=&\  \widetilde{R}_{(\ell,t)}+\sum_{i=1}^{N}\widetilde{D}^{i\top}_{(\ell,t)}\Phi_{(\ell,t)} \widetilde{D}^{i}_{(\ell,t)}, \quad &\ 
		\widetilde{L}_{(\ell,t)}=&\   \widetilde{S}_{(\ell,t)}+\Pi_{(\ell,t)}\widetilde{B}_{(\ell,t)}+\sum_{i=1}^{N}\widetilde{C}^{i\top}_{(\ell,t)}\Phi_{(\ell,t)} \widetilde{D}^{i}_{(\ell,t)}.
	\end{aligned}\right.
\end{equation}

The following lemma establishes the existence and uniqueness of solutions to \eqref{equ: Ricca one of Pro CMF}-\eqref{equ: Lya of Pro CMF}.
\begin{Lemma}\label{lem: wellpose of Ricc}
	Under assumptions \ref{ass:coef in com dyna}, \ref{ass: coe in com cost}, and \ref{ass: stand cod}, the Riccati equations \eqref{equ: Ricca one of Pro CMF} and \eqref{equ: Ricca two of Pro CMF} admit unique positive semidefinite solutions $\Phi_{(\ell,\cdot)}$ and $\Pi_{(\ell,\cdot)}$ such that $\Sigma_{(\ell,\cdot)}$ and $\widetilde{\Sigma}_{(\ell,\cdot)}$ are uniformly positive definite for every $\ell\in\mathbb{D}$. The BSDE \eqref{equ: Lya of Pro CMF} admits a unique solution $(\varphi,\rho)\in\mathcal{S}^{2}_{\mathbb{F}^{\theta}}([0,T];\mathbb{R}^{n})\times \mathcal{M}^{2}_{\mathbb{F}^{\theta}}([0,T];\mathbb{R}^{n})$.
\end{Lemma}
\begin{proof}
	The proof is provided in \ref{app: proof wellpose Ricc}.
\end{proof}

We now state the main result for \textbf{Problem CMF} under the positive definiteness assumption \ref{ass: stand cod}.
\begin{Theorem}\label{the: main result definite}
	Suppose that assumptions \ref{ass:coef in com dyna}, \ref{ass: coe in com cost}, and \ref{ass: stand cod} hold. Consider the following closed-loop system:
	\begin{equation}\label{equ: clos-loo sys CMF}
		\left\{\begin{aligned}
			\mathrm{d}X^{*}_t=&\ \bigl[\bigl(A-B \Sigma^{-1}L^{\top}\bigr)\bigl(X^{*}_t-\mathbb{E}^{\theta}_t\bigl[X^{*}_t\bigr]\bigr)+\bigl(\widetilde{A}-\widetilde{B} \widetilde{\Sigma}^{-1}\widetilde{L}^{\top}\bigr)\mathbb{E}^{\theta}_t\bigl[X^{*}_t\bigr]+\alpha-\widetilde{B}\widetilde{\Sigma}^{-1}\hbar\bigr]\mathrm{d}t\\
			&\  +\sum_{i=1}^{N}\Bigl[\bigl(C^{i}-D^{i} \Sigma^{-1}L^{\top}\bigr)\bigl(X^{*}_t-\mathbb{E}^{\theta}_t\bigl[X^{*}_t\bigr]\bigr)+\bigl(\widetilde{C}^{i}-\widetilde{D}^{i} \widetilde{\Sigma}^{-1}\widetilde{L}^{\top}\bigr)\mathbb{E}^{\theta}_t\bigl[X^{*}_t\bigr]+\beta^{i}-\widetilde{D}^{i}\widetilde{\Sigma}^{-1}\hbar\Bigr]\mathrm{d}W^i_t,\\
			X^{*}_0=&\ \xi,
		\end{aligned}\right.
	\end{equation}where \(			\hbar_{t}=\widetilde{B}_{(\theta_t,t)}^{\top}\varphi_{t}+\widetilde{r}_{t}+\sum_{i=1}^{N}\widetilde{D}^{i\top}_{(\theta_t,t)}\Phi_{(\theta_t,t)}\beta^{i}_{t}.\)	Define
	\begin{equation}\label{equ: u CMF}
		\begin{aligned}
			u^{*}_{t}=&\ -\Sigma^{-1}L^{\top}\bigl(X^{*}_t-\mathbb{E}^{\theta}_t\bigl[X^{*}_t\bigr]\bigr)-\widetilde{\Sigma}^{-1}\bigl(\widetilde{L}^{\top}\mathbb{E}^{\theta}_t\bigl[X^{*}_t\bigr]+\hbar_{t}\bigr), 
		\end{aligned}
	\end{equation}and
	\begin{equation}\label{equ: p q k CMF}
		\left\{\begin{aligned}
			p_t=&\  \Phi\bigl(X^{*}_t-\mathbb{E}^{\theta}_t\bigl[X^{*}_t\bigr]\bigr)+\Pi\mathbb{E}^{\theta}_t\bigl[X^{*}_t\bigr]+\varphi_{t},\\
			\zeta^{i}_t=&\ \Phi\bigl[\bigl(C^{i}-D^{i}\Sigma^{-1}L^{\top}\bigr)\bigl(X^{*}_t-\mathbb{E}^{\theta}_t\bigl[X^{*}_t\bigr]\bigr)+\bigl(\widetilde{C}^{i}-\widetilde{D}^{i}\widetilde{\Sigma}^{-1}\widetilde{L}^{\top}\bigr)\mathbb{E}^{\theta}_t\bigl[X^{*}_t\bigr]-\widetilde{D}^{i}\widetilde{\Sigma}^{-1}\hbar_{t}+\beta^{i}_{t}\bigr],\\
			k^{\ell j}_t=&\ \bigl(\Phi_{(j,t)}-\Phi_{(\ell,t)}\bigr)\bigl(X^{*}_{t}-\mathbb{E}^{\theta}_{t}\bigl[X^{*}_{t}\bigr]\bigr) +\bigl(\Pi_{(j,t)}-\Pi_{(\ell,t)}\bigr)\mathbb{E}^{\theta}_{t}\bigl[X^{*}_{t}\bigr]+\rho^{\ell j}_{t}.
		\end{aligned}\right.
	\end{equation}
	Then $(X^{*},u^{*},p,\zeta,k)$ satisfies the FBSDE \eqref{equ: Ham for com}, and $(X^{*},u^{*})$ is the optimal pair for \textbf{Problem CMF}.
	Moreover, the corresponding optimal cost is
	\begin{equation}\label{equ: exp cost of CMF}
		\begin{aligned}
			\mathscr{V}\bigl(\xi,\ell_0\bigr)
			=&\  \frac{1}{2}\mathbb{E}\biggl\{\int_{0}^{T}\Bigl[\sum_{i=1}^{N}\left|\beta^{i}_{t}\right|^{2}_{\Phi_{(\theta_t,t)}}-\left|\hbar_{t}\right|^{2}_{\widetilde{\Sigma}^{-1}_{(\theta_t,t)}}+2\alpha^{\top}_{t}\varphi_{t}\Bigr]\mathrm{d}t +\left|\xi-\mathbb{E}\xi\right|^{2}_{\Phi_{(\ell_0,0)}}+\left|\mathbb{E}\xi\right|^{2}_{\Pi_{(\ell_0,0)}}+2\bigl(\mathbb{E}\xi\bigr)^{\top}\varphi_{0}\biggr\}.
		\end{aligned}
	\end{equation}
\end{Theorem}
\begin{proof}
	The proof is provided in \ref{app: proof main definite}.
\end{proof}

\section{Indefinite Coefficients Case}\label{sec: CMF of indef coe}

Even when assumption~\ref{ass: stand cod} fails, \textbf{Problem CMF} may remain well posed and admit an optimal control. We address this case using the equivalent cost functional method. The core of this method is to replace the original cost functional with an equivalent one whose coefficients satisfy assumption~\ref{ass: stand cod}. We first define equivalence between cost functionals for a given controlled system.

\begin{Definition}\label{def: equivalent cost}
	Consider a controlled system with admissible control set $\mathcal U$.
	Two cost functionals $J$ and $\bar J$ are said to be
	\emph{equivalent} if they induce the same ordering on $\mathcal U$.
	That is, for any $u^1,u^2\in\mathcal U$,
	$J(u^1)<J(u^2)$ if and only if $\bar J(u^1)<\bar J(u^2)$.
\end{Definition}

By Definition~\ref{def: equivalent cost}, equivalent cost functionals for the same controlled system have the same set of optimal controls. Thus, when the coefficients of \(J\) fail to satisfy assumption~\ref{ass: stand cod}, we may seek an equivalent cost functional \(\bar{J}\) whose coefficients satisfy it. To construct such a functional, introduce the following class of deterministic matrix-valued function families:
\begin{equation*}
	\begin{aligned}
		\Delta=\bigl\{H=\bigl(H_{(\ell,\cdot)}\bigr)_{\ell\in\mathbb D}\ \big|\;& H_{(\ell,\cdot)}:\bigl[0,T\bigr]\rightarrow\mathbb S^n,
		H_{(\ell,\cdot)}\text{ is deterministic and absolutely continuous},\\
		& \dot H_{(\ell,\cdot)}\in L^\infty\bigl(\bigl[0,T\bigr];\mathbb S^n\bigr)
		\text{ for every }\ell\in\mathbb D\bigr\}.
	\end{aligned}
\end{equation*}
For any $H,K\in\Delta$, the quantity $\frac{1}{2}\mathbb{E}\Bigl\{\left|\xi-\mathbb{E}\xi\right|^{2}_{H_{(\ell_0,0)}}+\left|\mathbb{E}\xi\right|^{2}_{K_{(\ell_0,0)}}+\int_{0}^{T}\sum_{i=1}^{N}\left|\beta^{i}_{t}\right|^{2}_{H_{(\theta_t,t)}}\mathrm{d}t\Bigr\}$ is finite and independent of the control $u$. Subtracting this constant from the original cost functional \eqref{equ: com cos} therefore yields the equivalent cost functional
\begin{equation}\label{equ: J HK and J}
	\begin{aligned}
		J^{HK}\bigl(\xi,\ell_0; u\bigr)
		=&\ 	J\bigl(\xi,\ell_0; u\bigr)-\frac{1}{2}\mathbb{E}\biggl\{\left|\xi-\mathbb{E}\xi\right|^{2}_{H_{(\ell_0,0)}}+\left|\mathbb{E}\xi\right|^{2}_{K_{(\ell_0,0)}}+\int_{0}^{T}\sum_{i=1}^{N}\left|\beta^{i}_{t}\right|^{2}_{H_{(\theta_t,t)}}\mathrm{d}t\biggr\}.
	\end{aligned}
\end{equation}
By applying It\^o's formula to $\langle X^{u}_{t}-\mathbb{E}^{\theta}_t[X^{u}_{t}],H_{(\theta_t,t)}(X^{u}_{t}-\mathbb{E}^{\theta}_t[X^{u}_{t}]) \rangle+\langle \mathbb{E}^{\theta}_t[X^{u}_{t}],K_{(\theta_t,t)}\mathbb{E}^{\theta}_t[X^{u}_{t}] \rangle$, we can reformulate \eqref{equ: J HK and J} as
\begin{small}
	\begin{equation}\label{equ: com cos HK}
		\begin{aligned}
			J^{HK}\bigl(\xi,\ell_0; u\bigr)
			= &\  \frac{1}{2} \mathbb{E}\biggl\{\int _ { 0 } ^ { T } \biggl[\left\langle\left[\begin{array}{l}
				X^{u}_{t} \\
				u_{t}
			\end{array}\right],\left[\begin{array}{cc}
				Q^{HK}_{(\theta_t,t)} & S^{HK}_{(\theta_t,t)} \\
				\bigl(S^{HK}_{(\theta_t,t)}\bigr)^{\top} & R^{HK}_{(\theta_t,t)}
			\end{array}\right]\left[\begin{array}{l}
				X^{u}_{t} \\
				u_{t}
			\end{array}\right]\right\rangle+2\left\langle\left[\begin{array}{l}
				X^{u}_{t} \\
				u_{t}
			\end{array}\right],\left[\begin{array}{l}
				q^{HK}_{t} \\
				r^{HK}_{t}
			\end{array}\right]\right\rangle\\
			&\ +\left\langle\left[\begin{array}{l}
				\mathbb{E}^{\theta}_t\bigl[X^{u}_{t}\bigr] \\
				\mathbb{E}^{\theta}_t\bigl[u_{t}\bigr]
			\end{array}\right],\left[\begin{array}{cc}
				\widehat{Q}^{HK}_{(\theta_t,t)} & \widehat{S}^{HK}_{(\theta_t,t)} \\
				\bigl(\widehat{S}^{HK}_{(\theta_t,t)}\bigr)^{\top} & \widehat{R}^{HK}_{(\theta_t,t)}
			\end{array}\right]\left[\begin{array}{c}
				\mathbb{E}^{\theta}_t\bigl[X^{u}_{t}\bigr] \\
				\mathbb{E}^{\theta}_t\bigl[u_{t}\bigr]
			\end{array}\right]\right\rangle +2\left\langle\left[\begin{array}{l}
				\mathbb{E}^{\theta}_t\bigl[X^{u}_{t}\bigr] \\
				\mathbb{E}^{\theta}_t\bigl[u_{t}\bigr]
			\end{array}\right],\left[\begin{array}{l}
				\widehat{q}^{HK}_{t} \\
				\widehat{r}^{HK}_{t}
			\end{array}\right]\right\rangle\biggr]\mathrm{d} t\\
			&\  +\langle X^{u}_{T}, G^{HK}_{\theta_T} X^{u}_{T}\rangle+2\langle g, X^{u}_{T}\rangle+\langle\mathbb{E}^{\theta}_T\bigl[X^{u}_{T}\bigr], \widehat{G}^{HK}_{\theta_T} \mathbb{E}^{\theta}_T\bigl[X^{u}_{T}\bigr]\rangle+2\langle \widehat{g}, \mathbb{E}^{\theta}_T\bigl[X^{u}_{T}\bigr]\rangle\biggr\},
		\end{aligned}
	\end{equation}
\end{small}
where
\begin{small}
	\begin{equation}\label{equ: coe of J HK}
		\begin{aligned}
			Q^{HK}_{(\ell,t)}=&\ Q_{(\ell,t)}+\dot{H}_{(\ell,t)}+A^{\top}_{(\ell,t)}H_{(\ell,t)}+H_{(\ell,t)}A_{(\ell,t)}+\sum_{i=1}^{N}C^{i\top}_{(\ell,t)}H_{(\ell,t)}C^{i}_{(\ell,t)}+\sum_{j\in\mathbb{D}}\lambda^{\ell j}H_{(j,t)},\\
			S^{HK}_{(\ell,t)}=&\ S_{(\ell,t)}+H_{(\ell,t)}B_{(\ell,t)}+\sum_{i=1}^{N}C^{i\top}_{(\ell,t)}H_{(\ell,t)}D^{i}_{(\ell,t)},\quad 
			R^{HK}_{(\ell,t)}= R_{(\ell,t)}+\sum_{i=1}^{N}D^{i\top}_{(\ell,t)}H_{(\ell,t)}D^{i}_{(\ell,t)},\\
			q^{HK}_{t}=&\  q_{t}+\sum_{i=1}^{N}C^{i\top}_{(\theta_t,t)}H_{(\theta_t,t)}\beta^{i}_{t},\quad r^{HK}_{t}= r_{t}+\sum_{i=1}^{N}D^{i\top}_{(\theta_t,t)}H_{(\theta_t,t)}\beta^{i}_{t},\\
			\widetilde{Q}^{HK}_{(\ell,t)}=&\ \widetilde{Q}_{(\ell,t)}+\dot{K}_{(\ell,t)}+\widetilde{A}^{\top}_{(\ell,t)}K_{(\ell,t)}+K_{(\ell,t)}\widetilde{A}_{(\ell,t)}+\sum_{i=1}^{N}\widetilde{C}^{i\top}_{(\ell,t)}H_{(\ell,t)}\widetilde{C}^{i}_{(\ell,t)}+\sum_{j\in\mathbb{D}}\lambda^{\ell j}K_{(j,t)},\\
			\widetilde{S}^{HK}_{(\ell,t)}=&\ \widetilde{S}_{(\ell,t)}+K_{(\ell,t)}\widetilde{B}_{(\ell,t)}+\sum_{i=1}^{N}\widetilde{C}^{i\top}_{(\ell,t)}H_{(\ell,t)}\widetilde{D}^{i}_{(\ell,t)},\quad 
			\widetilde{R}^{HK}_{(\ell,t)}= \widetilde{R}_{(\ell,t)}+\sum_{i=1}^{N}\widetilde{D}^{i\top}_{(\ell,t)}H_{(\ell,t)}\widetilde{D}^{i}_{(\ell,t)},\\
			\widetilde{q}^{HK}_{t}=&\  \widetilde{q}_{t}+K_{(\theta_t,t)}\alpha_{t}+\sum_{i=1}^{N}\widetilde{C}^{i\top}_{(\theta_t,t)}H_{(\theta_t,t)}\beta^{i}_{t},\quad \widetilde{r}^{HK}_{t}= \widetilde{r}_{t}+\sum_{i=1}^{N}\widetilde{D}^{i\top}_{(\theta_t,t)}H_{(\theta_t,t)}\beta^{i}_{t}, \quad \widetilde{G}^{HK}_{\ell}= \widetilde{G}_{\ell}-K_{(\ell,T)},\\
			\widehat{Q}^{HK}_{(\ell,t)}=&\ \widetilde{Q}^{HK}_{(\ell,t)}-Q^{HK}_{(\ell,t)}, \quad	\widehat{S}^{HK}_{(\ell,t)}= \widetilde{S}^{HK}_{(\ell,t)}-S^{HK}_{(\ell,t)}, \quad	\widehat{R}^{HK}_{(\ell,t)}= \widetilde{R}^{HK}_{(\ell,t)}-R^{HK}_{(\ell,t)},\quad \widehat{q}^{HK}_{t}=\widetilde{q}^{HK}_{t}-q^{HK}_{t},\\
			\widehat{r}^{HK}_{t}=&\  \widetilde{r}^{HK}_{t}-r^{HK}_{t},\quad G^{HK}_{\ell}= G_{\ell}-H_{(\ell,T)},\quad \widehat{G}^{HK}_{\ell}=\widetilde{G}^{HK}_{\ell}-G^{HK}_{\ell}.\\
		\end{aligned}
	\end{equation}
\end{small}
The transformed coefficients in \eqref{equ: coe of J HK} may satisfy assumption~\ref{ass: stand cod} even when the original coefficients do not. A suitable choice of $H, K\in\Delta$ can therefore transform the indefinite \textbf{Problem CMF} into an LQ problem satisfying this assumption.

For given $H,K$, the state equation \eqref{equ: com sta} and cost functional \eqref{equ: com cos HK} define \textbf{Problem CMF$^{HK}$}. This problem reduces to \textbf{Problem CMF} when $H=K=0$. Its associated Hamiltonian system, coupled Riccati equations, and linear BSDE are given by
\begin{equation}\label{equ: Ham for HK}
	\left\{\begin{aligned}
		\mathrm{d}X^{HK}_t=&\ \bigl[AX^{HK}_t+\widehat{A}\mathbb{E}^{\theta}_t\bigl[X^{HK}_t\bigr]+Bu^{HK}_{t}+\widehat{B}\mathbb{E}^{\theta}_t\bigl[u^{HK}_{t}\bigr]+\alpha_{t}\bigr]\mathrm{d}t\\
		&\  +\sum_{i=1}^{N}\Bigl[C^{i}X^{HK}_t+\widehat{C}^{i}\mathbb{E}^{\theta}_t\bigl[X^{HK}_t\bigr]+D^{i}u^{HK}_{t}+\widehat{D}^{i}\mathbb{E}^{\theta}_t\bigl[u^{HK}_{t}\bigr]+\beta^{i}_{t}\Bigr]\mathrm{d}W^i_t,\\
		\mathrm{d}p^{HK}_t=&\   -\Bigl[A^{\top}p^{HK}_t+\widehat{A}^{\top}\mathbb{E}^{\theta}_t\bigl[p^{HK}_t\bigr]+\sum_{i=1}^{N}C^{i\top}\zeta^{HK,i}_t+\sum_{i=1}^{N}\widehat{C}^{i\top}\mathbb{E}^{\theta}_t\bigl[\zeta^{HK,i}_t\bigr]+Q^{HK}X^{HK}_t +\widehat{Q}^{HK}\mathbb{E}^{\theta}_t\bigl[X^{HK}_t\bigr]\\
		&\  +S^{HK}u^{HK}_{t} +\widehat{S}^{HK}\mathbb{E}^{\theta}_t\bigl[u^{HK}_{t}\bigr]+\widetilde{q}^{HK}_{t}\Bigr]\mathrm{d}t +\sum_{i=1}^{N}\zeta^{HK,i}_t\mathrm{d}W^{i}_t+k^{HK}_t\bullet\mathrm{d}M_t,\\
		X^{HK}_0=&\ \xi,\quad p^{HK}_T=G^{HK}_{\theta_T}X^{HK}_T+\widehat{G}^{HK}_{\theta_T}\mathbb{E}^{\theta}_T\bigl[X^{HK}_T\bigr]+\widetilde{g},\\
		&\ B^{\top}p^{HK}_t+\widehat{B}^{\top}\mathbb{E}^{\theta}_t\bigl[p^{HK}_t\bigr]+\sum_{i=1}^{N}D^{i\top}\zeta^{HK,i}_t+\sum_{i=1}^{N}\widehat{D}^{i\top}\mathbb{E}^{\theta}_t\bigl[\zeta^{HK,i}_t\bigr]\\
		&\ +\bigl(S^{HK}\bigr)^{\top}X^{HK}_t+\bigl(\widehat{S}^{HK}\bigr)^{\top}\mathbb{E}^{\theta}_t\bigl[X^{HK}_t\bigr]+R^{HK}u^{HK}_{t}+\widehat{R}^{HK}\mathbb{E}^{\theta}_t\bigl[u^{HK}_{t}\bigr]+\widetilde{r}^{HK}_{t}=0,
	\end{aligned}\right.
\end{equation}
\begin{equation}\label{equ: Ricca one of Pro CMF HK}
	\left\{\begin{aligned}
		\dot{\Phi}^{HK}_{(\ell,t)}=&\ -\Bigl[A^{\top}_{(\ell,t)}\Phi^{HK}_{(\ell,t)}+\Phi^{HK}_{(\ell,t)} A_{(\ell,t)}+\sum_{i=1}^{N}C^{i\top}_{(\ell,t)}\Phi^{HK}_{(\ell,t)} C^{i}_{(\ell,t)}+Q^{HK}_{(\ell,t)}\\
		&\ -L^{HK}_{(\ell,t)}\bigl(\Sigma^{HK}_{(\ell,t)}\bigr)^{-1}\bigl(L^{HK}_{(\ell,t)}\bigr)^{\top}+\sum_{j\in\mathbb{D}}\lambda^{\ell j}\Phi^{HK}_{(j,t)}\Bigr],\\
		\Phi^{HK}_{(\ell,T)}=&\ G^{HK}_{\ell},
	\end{aligned}\right.
\end{equation}	

\begin{equation}\label{equ: Ricca two of Pro CMF HK}
	\left\{\begin{aligned}
		\dot{\Pi}^{HK}_{(\ell,t)}=&\ -\Bigl[\widetilde{A}^{\top}_{(\ell,t)}\Pi^{HK}_{(\ell,t)}+\Pi^{HK}_{(\ell,t)} \widetilde{A}_{(\ell,t)}+\sum_{i=1}^{N}\widetilde{C}^{i\top}_{(\ell,t)}\Phi^{HK}_{(\ell,t)} \widetilde{C}^{i}_{(\ell,t)}+\widetilde{Q}^{HK}_{(\ell,t)}\\
		&\ -\widetilde{L}^{HK}_{(\ell,t)}\bigl(\widetilde{\Sigma}^{HK}_{(\ell,t)}\bigr)^{-1}\bigl(\widetilde{L}^{HK}_{(\ell,t)}\bigr)^{\top}+\sum_{j\in\mathbb{D}}\lambda^{\ell j}\Pi^{HK}_{(j,t)}\Bigr],\\
		\Pi^{HK}_{(\ell,T)}=&\ \widetilde{G}^{HK}_{\ell},
	\end{aligned}\right.
\end{equation}
and 
\begin{equation}\label{equ: Lya of Pro CMF HK}
	\left\{\begin{aligned}
		\mathrm{d}\varphi^{HK}_{t}=&\ -\Bigl[\bigl(\widetilde{A}-\widetilde{B}\bigl(\widetilde{\Sigma}^{HK}\bigr)^{-1}\bigl(\widetilde{L}^{HK}\bigr)^{\top}\bigr)^{\top}\varphi^{HK}_{t}+\sum_{i=1}^{N}\bigl(\widetilde{C}^{i}-\widetilde{D}^{i}\bigl(\widetilde{\Sigma}^{HK}\bigr)^{-1}\bigl(\widetilde{L}^{HK}\bigr)^{\top}\bigr)^{\top}\Phi^{HK}\beta^{i}_{t}\\
		&\  -\widetilde{L}^{HK}\bigl(\widetilde{\Sigma}^{HK}\bigr)^{-1} \widetilde{r}^{HK}_{t}+\widetilde{q}^{HK}_{t}+\Pi^{HK}\alpha_{t}\Bigr]\mathrm{d}t+\rho^{HK}_t\bullet\mathrm{d}M_t,\\
		\varphi^{HK}_{T}=&\ \widetilde{g},
	\end{aligned}\right.
\end{equation}
where 
\begin{equation*}
	\left\{\begin{aligned}
		\Sigma^{HK}_{(\ell,t)}=&\ R^{HK}_{(\ell,t)}+\sum_{i=1}^{N}D^{i\top}_{(\ell,t)}\Phi^{HK}_{(\ell,t)} D^{i}_{(\ell,t)}, \quad L^{HK}_{(\ell,t)}= S^{HK}_{(\ell,t)}+\Phi^{HK}_{(\ell,t)}B_{(\ell,t)}+\sum_{i=1}^{N}C^{i\top}_{(\ell,t)}\Phi^{HK}_{(\ell,t)} D^{i}_{(\ell,t)},\\
		\widetilde{\Sigma}^{HK}_{(\ell,t)}=&\ \widetilde{R}^{HK}_{(\ell,t)}+\sum_{i=1}^{N}\widetilde{D}^{i\top}_{(\ell,t)}\Phi^{HK}_{(\ell,t)} \widetilde{D}^{i}_{(\ell,t)},\quad 
		\widetilde{L}^{HK}_{(\ell,t)}= \widetilde{S}^{HK}_{(\ell,t)}+\Pi^{HK}_{(\ell,t)}\widetilde{B}_{(\ell,t)}+\sum_{i=1}^{N}\widetilde{C}^{i\top}_{(\ell,t)}\Phi^{HK}_{(\ell,t)} \widetilde{D}^{i}_{(\ell,t)}.
	\end{aligned}\right.
\end{equation*}

The following lemma relates the Hamiltonian systems of \textbf{Problem CMF$^{HK}$} and \textbf{Problem CMF}.
\begin{Lemma}\label{lem: Ham of indefi}
	Let $H,K\in\Delta$, and suppose that
	$(X^{HK},u^{HK},p^{HK},\zeta^{HK},k^{HK})$
	is a solution of the Hamiltonian system \eqref{equ: Ham for HK}.
	Define $(X,u,p,\zeta,k)$ by
	\begin{equation}\label{equ: Ham inve trans}
		\left\{\begin{aligned}
			X_t=&\ X^{HK}_t, \quad u_t=\ u^{HK}_{t},\quad 
			p_t=  p^{HK}_t+H\bigl(X^{HK}_t-\mathbb{E}^{\theta}_t\bigl[X^{HK}_t\bigr]\bigr)+K\mathbb{E}^{\theta}_t\bigl[X^{HK}_t\bigr],\\
			\zeta^i_t=&\ \zeta^{HK,i}_t+H\bigl(C^{i}X^{HK}_t+\widehat{C}^{i}\mathbb{E}^{\theta}_t\bigl[X^{HK}_t\bigr]+D^{i}u^{HK}_{t}+\widehat{D}^{i}\mathbb{E}^{\theta}_t\bigl[u^{HK}_{t}\bigr]+\beta^{i}_{t}\bigr),\\
			k^{\ell j}_t=&\ k^{HK,\ell j}_t+\bigl(H_{(j,t)}-H_{(\ell,t)}\bigr)\bigl(X^{HK}_t-\mathbb{E}^{\theta}_t\bigl[X^{HK}_t\bigr]\bigr)+\bigl(K_{(j,t)}-K_{(\ell,t)}\bigr)\mathbb{E}^{\theta}_t\bigl[X^{HK}_t\bigr].
		\end{aligned}\right.
	\end{equation}
	Then $(X,u,p,\zeta,k)$ solves \eqref{equ: Ham for com}. Moreover, for every $H,K\in\Delta$, the existence and uniqueness of solutions to  systems \eqref{equ: Ham for HK} and \eqref{equ: Ham for com} are equivalent.
\end{Lemma}
\begin{proof}
	Apply It\^o's formula to
	\(
	H_{(\theta_t,t)}
	(
	X^{HK}_t-\mathbb{E}^{\theta}_t[X^{HK}_t]
	)
	+K_{(\theta_t,t)}
	\mathbb{E}^{\theta}_t[X^{HK}_t].
	\)
	Combining this identity with the dynamics of $p^{HK}$
	and the coefficient relations in \eqref{equ: coe of J HK}
	yields the backward equation in \eqref{equ: Ham for com}.
	Direct substitution verifies the state equation,
	terminal condition, and maximum condition
	in \eqref{equ: Ham for com}.
	Thus, $(X,u,p,\zeta,k)$ is a solution of
	\eqref{equ: Ham for com}. Conversely, \eqref{equ: Ham inve trans} is an invertible transformation.
	The same calculation verifies that the  inverse transformation
	maps every solution of \eqref{equ: Ham for com}
	to a solution of \eqref{equ: Ham for HK}. Thus, \eqref{equ: Ham inve trans} defines a bijection between the two solution sets
	and preserves both existence and uniqueness.
\end{proof}
An analogous conclusion holds for the Riccati equations and the associated linear BSDE. The next lemma can be verified by direct substitution using \eqref{equ: coe of J HK}, so its proof is omitted.
\begin{Lemma}\label{lem: Ricca of indefi}
	Let $H,K\in\Delta$.
	Suppose that $\Phi^{HK}_{(\ell,\cdot)}$ and $\Pi^{HK}_{(\ell,\cdot)}$
	are solutions of
	\eqref{equ: Ricca one of Pro CMF HK} and
	\eqref{equ: Ricca two of Pro CMF HK}, respectively,
	and that $(\varphi^{HK},\rho^{HK})$ solves
	\eqref{equ: Lya of Pro CMF HK}.
	Define $\Phi$, $\Pi$, and $(\varphi,\rho)$ by
	\begin{equation*}
		\begin{aligned}
			\Phi_{(\ell,t)}
			= \Phi^{HK}_{(\ell,t)}+H_{(\ell,t)},\quad 
			\Pi_{(\ell,t)}
			= \Pi^{HK}_{(\ell,t)}+K_{(\ell,t)},\quad 
			\bigl(\varphi_t,\rho_t\bigr)
			= \bigl(\varphi^{HK}_t,\rho^{HK}_t\bigr).
		\end{aligned}
	\end{equation*}
	Then $\Phi_{(\ell,\cdot)}$ and $\Pi_{(\ell,\cdot)}$ are solutions of
	\eqref{equ: Ricca one of Pro CMF} and
	\eqref{equ: Ricca two of Pro CMF}, respectively,
	and $(\varphi,\rho)$ solves \eqref{equ: Lya of Pro CMF}.
	Moreover, for every $H,K\in\Delta$, the existence and
	uniqueness of solutions to
	\eqref{equ: Ricca one of Pro CMF HK}-\eqref{equ: Lya of Pro CMF HK}
	are equivalent to those for
	\eqref{equ: Ricca one of Pro CMF}-\eqref{equ: Lya of Pro CMF}.
\end{Lemma}

We now establish the main result for the indefinite coefficients case.

\begin{Theorem}\label{the: H bar K bar}
	Suppose that assumptions \ref{ass:coef in com dyna}
	and \ref{ass: coe in com cost} hold.
	Assume that there exist $\bar{H},\bar{K}\in\Delta$
	such that \textbf{Problem CMF$^{\bar{H}\bar{K}}$}
	satisfies assumption~\ref{ass: stand cod}. Then the Hamiltonian system \eqref{equ: Ham for com} admits a unique solution $(X^*,u^*,p,\zeta,k)$. The Riccati equations
	\eqref{equ: Ricca one of Pro CMF}-\eqref{equ: Ricca two of Pro CMF} and the regime-switching BSDE \eqref{equ: Lya of Pro CMF} admit unique solutions $\Phi_{(\ell,\cdot)}$, $\Pi_{(\ell,\cdot)}$, and $(\varphi,\rho)$.
	Moreover, $\Sigma_{(\ell,\cdot)}$ and
	$\widetilde{\Sigma}_{(\ell,\cdot)}$ are uniformly
	positive definite for every $\ell\in\mathbb{D}$.
	
	Let $X^*$ solve \eqref{equ: clos-loo sys CMF},
	and define $(u^*,p,\zeta,k)$ by
	\eqref{equ: u CMF}-\eqref{equ: p q k CMF}.
	Then $(X^*,u^*,p,\zeta,k)$ satisfies
	\eqref{equ: Ham for com},
	and $(X^*,u^*)$ is the unique optimal pair
	of \textbf{Problem CMF}.
	The corresponding optimal cost is given by
	\eqref{equ: exp cost of CMF}.
\end{Theorem}
\begin{proof}
	The proof is provided in \ref{app: proof H bar K bar}.
\end{proof}

\begin{Lemma}\label{lem: exis of HK}
	Suppose that assumptions~\ref{ass:coef in com dyna}
	and~\ref{ass: coe in com cost} hold.
	Assume that the Riccati equations
	\eqref{equ: Ricca one of Pro CMF} and
	\eqref{equ: Ricca two of Pro CMF}
	admit solutions
	$\Phi_{(\ell,\cdot)}$ and $\Pi_{(\ell,\cdot)}$, respectively, such that
	$\Sigma_{(\ell,\cdot)}$ and $\widetilde{\Sigma}_{(\ell,\cdot)}$ are uniformly
	positive definite. Then choosing
	$\bar H=\Phi$ and $\bar K=\Pi$ ensures that
	\textbf{Problem CMF$^{\bar H\bar K}$}
	satisfies assumption~\ref{ass: stand cod}.
\end{Lemma}
\begin{proof}
	Set $\bar H=\Phi$ and $\bar K=\Pi$.
	By \eqref{equ: coe of J HK} and \eqref{equ: nota Sigma}, we have
	\begin{equation*}
		\begin{aligned}
			R^{\Phi\Pi}_{(\ell,t)}
			=\Sigma_{(\ell,t)},\quad  S^{\Phi\Pi}_{(\ell,t)}
			=L_{(\ell,t)}, \quad 
			\widetilde R^{\Phi\Pi}_{(\ell,t)}
			=\widetilde\Sigma_{(\ell,t)},\quad  
			\widetilde S^{\Phi\Pi}_{(\ell,t)}
			=\widetilde L_{(\ell,t)}.
		\end{aligned}
	\end{equation*}  Furthermore, the Riccati equations \eqref{equ: Ricca one of Pro CMF} and
	\eqref{equ: Ricca two of Pro CMF} yield
	\begin{equation*}
		\begin{aligned}
			Q^{\Phi\Pi}_{(\ell,t)}
			=
			L_{(\ell,t)}
			\Sigma_{(\ell,t)}^{-1}
			L_{(\ell,t)}^\top, \quad  
			\widetilde Q^{\Phi\Pi}_{(\ell,t)}
			=
			\widetilde L_{(\ell,t)}
			\widetilde\Sigma_{(\ell,t)}^{-1}
			\widetilde L_{(\ell,t)}^\top.
		\end{aligned}
	\end{equation*}
	Consequently, $Q^{\Phi\Pi}_{(\ell,t)}
	-
	S^{\Phi\Pi}_{(\ell,t)}
	(R^{\Phi\Pi}_{(\ell,t)})^{-1}
	(S^{\Phi\Pi}_{(\ell,t)})^\top=0$, $
	\widetilde Q^{\Phi\Pi}_{(\ell,t)}
	-
	\widetilde S^{\Phi\Pi}_{(\ell,t)}
	(\widetilde R^{\Phi\Pi}_{(\ell,t)})^{-1}
	(\widetilde S^{\Phi\Pi}_{(\ell,t)})^\top
	=0.$
	
	By the assumed uniform positive definiteness,
	there exists a constant $\delta>0$ such that $R^{\Phi\Pi}_{(\ell,t)}
	\geq\delta\mathbf{I}_{m\times m}$, $
	\widetilde R^{\Phi\Pi}_{(\ell,t)}
	\geq\delta\mathbf{I}_{m\times m}$. Finally, the terminal conditions give
	$
	G^{\Phi\Pi}_\ell
	=G_\ell-\Phi_{(\ell,T)}=0$, 
	$\widetilde G^{\Phi\Pi}_\ell
	=\widetilde G_\ell-\Pi_{(\ell,T)}=0$.
	Hence, the transformed coefficients satisfy all the conditions in
	assumption~\ref{ass: stand cod}.
\end{proof}

\section{Examples}\label{sec:example}

This section illustrates the equivalent cost functional method through two
examples. In the first example,  the coefficients violate assumption~\ref{ass: stand cod}. However, the equivalent cost functional method yields an optimal control in feedback form and an
explicit expression for the optimal cost. In the second example, we identify a fully coupled CMF-FBSDE  with the Hamiltonian system of an indefinite LQ control problem. By applying Theorem \ref{the: H bar K bar}, we establish the  unique solvability
of the  FBSDE  even though the monotonicity conditions examined below are not satisfied.

\medskip
\noindent\textbf{Example 1.}
Let $\theta$ be a Markov chain with state space
$\mathbb D=\{1,2\}$, initial state $\theta_0=\ell_0$, and transition intensity  \(\lambda^{12}=\lambda^{21}=1\). Let $W$ be a one-dimensional Brownian motion independent of $\theta$.
Consider the scalar controlled system
\[
\begin{aligned}
	\mathrm dX_t={}
	\bigl[A_{\theta_t}\bigl(X_t-\mathbb E_t^\theta\bigl[X_t\bigr]\bigr)+u_t\bigr]
	\,\mathrm dt+\bigl[2\bigl(u_t-\mathbb E_t^\theta\bigl[u_t\bigr]\bigr)
	+\mathbb E_t^\theta\bigl[u_t\bigr]\bigr]\,\mathrm dW_t,\quad X_0={}x,
\end{aligned}
\]
where $x\in\mathbb R$ is deterministic, $A_1=0$, and $A_2=\frac{1}{3}$.
The cost functional is
\[
\begin{aligned}
	J\bigl(x,\ell_0;u\bigr)=\frac12\mathbb E\biggl\{
	&G_{\theta_T}X_T^2
	+\widehat G_{\theta_T}\bigl(\mathbb E_T^\theta\bigl[X_T\bigr]\bigr)^2+\int_0^T\Bigl[
	Q_{\theta_t}X_t^2-u_t^2
	+\widehat Q_{\theta_t}\bigl(\mathbb E_t^\theta\bigl[X_t\bigr]\bigr)^2
	+\frac12\bigl(\mathbb E_t^\theta\bigl[u_t\bigr]\bigr)^2
	\Bigr]\,\mathrm dt\biggr\},
\end{aligned}
\]
with coefficients $G_{1}=1$, $G_{2}=\frac{3}{2}$, $\widehat{G}_{1}=-\frac{1}{2}$, $\widehat{G}_{2}=-\frac{3}{4}$, $Q_{1}=-\frac{1}{6}$,  $Q_{2}=-\frac{1}{20}$, $\widehat{Q}_{1}=\frac{5}{12}$, and $\widehat{Q}_{2}=\frac{69}{80}$.  For both regimes, $R_\ell=-1$ and $\widetilde R_\ell=-1/2$.
Consequently, the original problem does not satisfy
assumption~\ref{ass: stand cod}.

The Riccati equations \eqref{equ: Ricca one of Pro CMF} and
\eqref{equ: Ricca two of Pro CMF} reduce to
\begin{equation*}
	\left\{\begin{aligned}
		\dot\Phi_{(1,t)}
		&=-\biggl[-\frac16
		-\frac{\Phi_{(1,t)}^2}{4\Phi_{(1,t)}-1}
		-\Phi_{(1,t)}+\Phi_{(2,t)}\biggr],\quad &\Phi_{(1,T)}=&\ 1,\\
		\dot\Phi_{(2,t)}
		&=-\biggl[\Phi_{(1,t)}-\frac13\Phi_{(2,t)}-\frac1{20}
		-\frac{\Phi_{(2,t)}^2}{4\Phi_{(2,t)}-1}\biggr],\quad & \Phi_{(2,T)}=&\ \frac32, \\
	\end{aligned}\right.
\end{equation*}and
\begin{equation*}
	\left\{\begin{aligned}
		\dot{\Pi}_{(1,t)} &= -\biggl[ \frac{1}{4} - \frac{\Pi^2_{(1,t)}}{\Phi_{(1,t)}-\frac12} - \Pi_{(1,t)} + \Pi_{(2,t)} \biggr],\quad & \Pi_{(1,T)} &= \frac12,  \\
		\dot{\Pi}_{(2,t)} &= -\biggl[ \frac{13}{16} - \frac{\Pi^2_{(2,t)}}{\Phi_{(2,t)}-\frac12} + \Pi_{(1,t)} - \Pi_{(2,t)} \biggr],\quad & \Pi_{(2,T)} &= \frac34.
	\end{aligned}\right.
\end{equation*}
We can verify  that these equations admit the constant solutions 
\(
\Phi_{(1,t)}=1,\Phi_{(2,t)}=\frac32,
\Pi_{(1,t)}=\frac12,\Pi_{(2,t)}=\frac34.
\)
The corresponding  matrices satisfy
\(
\Sigma_{(1,t)}=3, \Sigma_{(2,t)}=5,
\widetilde\Sigma_{(1,t)}=\frac12,
\widetilde\Sigma_{(2,t)}=1.
\)
Thus, both $\Sigma$ and $\widetilde\Sigma$ are uniformly positive definite.
By Lemma~\ref{lem: exis of HK}, the choice $\bar H=\Phi$ and
$\bar K=\Pi$ yields a transformed problem satisfying
assumption~\ref{ass: stand cod}.
Theorem~\ref{the: H bar K bar} therefore gives the unique optimal control
in feedback form.

Define
\(
a_1=\frac13, a_2=\frac3{10},
b_1=1, b_2=\frac34.
\) It follows from  Theorem~\ref{the: H bar K bar} that the optimal control and associated cost are  given by 
\[u_t^*={}-a_{\theta_t}\bigl(X_t^*-\mathbb E_t^\theta\bigl[X_t^*\bigr]\bigr)-b_{\theta_t}\mathbb E_t^\theta\bigl[X_t^*\bigr],\] and 
\[
\mathscr V\bigl(x,\ell_0\bigr)
=\frac12\Pi_{(\ell_0,0)}x^2
=\begin{cases}
	\displaystyle \frac14x^2,\quad &\ell_0=1,\\[5pt]
	\displaystyle \frac38x^2,\quad &\ell_0=2.
\end{cases}
\]
%

\medskip
\noindent\textbf{Example 2.}
Our theoretical result provides an effective tool for establishing the well-posedness of fully coupled
CMF-FBSDEs with regime switching.  Establishing global well-posedness for this class of equations is challenging. To the best of our knowledge, global well-posedness of such systems over finite time horizons has been studied only in
\cite{Rol贸n Guti茅rrez markovian 2024,HuangLiXuZheng2026,Chen Shi Wang Wu 2026,WangChakraborty2026}.
However, the results in these references do not establish
the solvability of \eqref{eq:scalar-fbsde}.
Specifically, as verified below, the
monotonicity conditions in \cite{Rol贸n Guti茅rrez markovian 2024} and  \cite{HuangLiXuZheng2026}
are not satisfied by \eqref{eq:scalar-fbsde}. 
The existence and uniqueness assertion in
\cite[Theorem~4.10]{WangChakraborty2026} is established by applying
\cite[Theorems~3.1 and~3.2]{Rol贸n Guti茅rrez markovian 2024}.
Consequently, the failure of the required monotonicity conditions
also precludes the application of Wang and Chakraborty's
result to \eqref{eq:scalar-fbsde}.
The analysis in \cite{Chen Shi Wang Wu 2026} is motivated by the
Hamiltonian system of a backward mean-field game and concerns an
FBSDE with coupling in the initial  condition and a
prescribed terminal value for the backward equation.
This boundary structure differs from that of \eqref{eq:scalar-fbsde},
in which the terminal condition of the backward equation depends on the solution of the forward SDE. Fortunately,  by identifying  \eqref{eq:scalar-fbsde} with the Hamiltonian systems of indefinite LQ control problems and using Theorem \ref{the: H bar K bar}, we prove the   solvability of \eqref{eq:scalar-fbsde}.

Let $T>0$, and let $W$ be a one-dimensional Brownian motion independent
of a Markov chain $\theta$ with state space $\mathbb D=\{1,2\}$,
initial state $\theta_0=\ell_0$, and transition intensities \(\lambda^{12}=1\) and \(\lambda^{21}=2\). Consider the scalar CMF-FBSDE
\begin{equation}
	\label{eq:scalar-fbsde}
	\left\{\begin{aligned}
		\mathrm{d}X_t={}&\bigl(X_t+\mathbb{E}^{\theta}_t\bigl[X_t\bigr]
		+\gamma_{\theta_t}p_t+\gamma_{\theta_t}\zeta_t\bigr)\,\mathrm{d}t+\bigl(-X_t
		+\gamma_{\theta_t}p_t+\gamma_{\theta_t}\zeta_t\bigr)\,\mathrm{d}W_t,\quad X_0={}x,\\
		\mathrm{d}p_t={}&\bigl(X_t+\mathbb{E}^{\theta}_t\bigl[X_t\bigr]-p_t-\mathbb{E}^{\theta}_t\bigl[p_t\bigr]+\zeta_t\bigr)\,\mathrm{d}t+\zeta_t\,\mathrm{d}W_t+k_t\bullet \mathrm{d}M_t,\qquad p_T=G_{\theta_T}X_T+\mathbb{E}^{\theta}_T\bigl[X_T\bigr],
	\end{aligned}\right.
\end{equation}
where $x\in\mathbb R$, $\gamma_1=1$, $\gamma_2=-\frac{1}{2}$,
$G_1=2$, and $G_2=-1$.

Using the sign convention  in
\cite{Rol贸n Guti茅rrez markovian 2024}, we write 
\(h(x,\mu,\ell)=G_\ell x+\int_{\mathbb R}\xi\,\mu(\mathrm d\xi)\).
For any fixed $\mu$ and regime $\ell=2$, the terminal value satisfies
\[
\bigl(h\bigl(x_1,\mu,2\bigr)-h\bigl(x_2,\mu,2\bigr)\bigr)\bigl(x_1-x_2\bigr)
=-|x_1-x_2|^2.
\]
This violates condition (H1) in
\cite{Rol贸n Guti茅rrez markovian 2024}, which requires
a lower bound of $K_h|x_1-x_2|^2$ for some $K_h>0$.



We next examine the monotonicity conditions underlying
Theorem~2.1 of \cite{HuangLiXuZheng2026}.
For each $\ell\in\mathbb D$, let $b_\ell$, $\sigma_\ell$,
and $g_\ell$ denote the drift,  diffusion
coefficient, and generator in \eqref{eq:scalar-fbsde},
respectively.
Set $\Gamma_\ell=(g_\ell,b_\ell,\sigma_\ell)^\top$.
For two square-integrable $\mathcal F_t$-measurable triples
$\Theta=(X,p,\zeta)^\top$ and
$\Theta'=(X',p',\zeta')^\top$, write
\(
\widehat\Theta=\mathbb E_t^\theta\bigl[\Theta\bigr],
\widehat\Theta'=\mathbb E_t^\theta\bigl[\Theta'\bigr],
\Delta\Theta=\Theta-\Theta'.
\)
A direct calculation, with the conditional expectation evaluated
at the respective triples, gives
\[
\begin{aligned}
	\mathbb E_t^\theta\bigl[
	\bigl\langle
	\Gamma_\ell\bigl(\Theta,\widehat\Theta\bigr)
	-\Gamma_\ell\bigl(\Theta',\widehat\Theta'\bigr),
	\Delta\Theta
	\bigr\rangle
	\bigr] =
	\mathbb E_t^\theta\bigl[\bigl(\Delta X\bigr)^2\bigr]
	+\bigl(\mathbb E_t^\theta\bigl[\Delta X\bigr]\bigr)^2
	+\gamma_\ell\mathbb E_t^\theta\bigl[
	\bigl(\Delta p+\Delta\zeta\bigr)^2
	\bigr].
\end{aligned}
\]
Taking $\Delta X=1$ and $\Delta p=\Delta\zeta=0$ makes
this expression equal to $2$.
This contradicts the nonpositivity required by the last
inequality in (2.4) of Assumption~(A2) in
\cite{HuangLiXuZheng2026}. Identifying \eqref{eq:scalar-fbsde} with the second block
of that framework does not resolve this issue.
Indeed, for $\ell=2$, taking $\Delta X=\Delta\zeta=0$
and $\Delta p=1$ gives $\gamma_2=-1/2$.
This contradicts the nonnegativity required by the last
inequality in (2.5) of the same assumption.
Consequently, Theorem~2.1 of \cite{HuangLiXuZheng2026}
cannot be applied  to 
\eqref{eq:scalar-fbsde}.

We now associate \eqref{eq:scalar-fbsde} with an indefinite control problem.
Consider the controlled state equation
\begin{equation*}\label{eq:scalar-state}
	\begin{aligned}
		\mathrm dX_t={}
		\bigl(A_{\theta_t}X_t+\mathbb E_t^\theta\bigl[X_t\bigr]+u_t\bigr)
		\,\mathrm dt+\bigl(X_t+u_t\bigr)\,\mathrm dW_t,\quad 
		X_0={}x,
	\end{aligned}
\end{equation*}
and the cost functional
\begin{equation*}\label{eq:scalar-cost}
	\begin{aligned}
		J\bigl(x,\ell_0;u\bigr)=\frac12\mathbb E\Bigl\{
		G_{\theta_T}X_T^2+\bigl(\mathbb E_T^\theta\bigl[X_T\bigr]\bigr)^2+\int_0^T\Bigl[
		Q_{\theta_t}X_t^2+2S_{\theta_t}X_tu_t+R_{\theta_t}u_t^2
		+\widehat Q_{\theta_t}\bigl(\mathbb E_t^\theta\bigl[X_t\bigr]\bigr)^2
		\Bigr]\,\mathrm dt\Bigr\},
	\end{aligned}
\end{equation*}
where  $A_{1}=A_{2}=3$, $Q_{1}=-5$, $Q_{2}=7$, $\widehat{Q}_{1}=\widehat{Q}_{2}=-1$, $S_{1}=-2$, $S_{2}=4$, $R_{1}=-1$, $R_{2}=2$, $G_{1}=2$, and $G_{2}=-1$.

The associated Hamiltonian system is
\begin{equation}\label{equ: exa 2 FBSDE}
	\left\{
	\begin{aligned}
		\mathrm dX_t={}&
		\bigl(A_{\theta_t}X_t+\mathbb E_t^\theta\bigl[X_t\bigr]+u_t\bigr)
		\,\mathrm dt+\bigl(X_t+u_t\bigr)\,\mathrm dW_t,\\
		\mathrm dp_t={}&
		-\bigl(A_{\theta_t}p_t+\mathbb E_t^\theta\bigl[p_t\bigr]+\zeta_t
		+Q_{\theta_t}X_t+\widehat Q_{\theta_t}\mathbb E_t^\theta\bigl[X_t\bigr]
		+S_{\theta_t}u_t\bigr)\,\mathrm dt
		+\zeta_t\,\mathrm dW_t+k_t\bullet\mathrm dM_t,\\
		X_0={}&x,\quad 
		p_T=G_{\theta_T}X_T+\mathbb E_T^\theta\bigl[X_T\bigr],\\
		&\hspace{-0.5cm}p_t+\zeta_t+S_{\theta_t}X_t+R_{\theta_t}u_t=0.
	\end{aligned}
	\right.
\end{equation}
Since $R_\ell\ne0$ in both regimes, the maximum condition gives
\(
u_t=-R_{\theta_t}^{-1}
(p_t+\zeta_t+S_{\theta_t}X_t).
\) Substituting this expression for $u_t$ into the state and adjoint equations
gives exactly \eqref{eq:scalar-fbsde}.
Conversely, every square-integrable adapted solution of
\eqref{eq:scalar-fbsde} determines a solution of
\eqref{equ: exa 2 FBSDE} through the same expression for $u_t$.

The original control problem violates assumption~\ref{ass: stand cod},
since $R_1=-1$ and $G_2=-1$.
Choose the constant auxiliary functions \(H_{(1,t)}=\frac32, H_{(2,t)}=-1,
K_{(1,t)}=\frac52, K_{(2,t)}=0\),  \(t\in\bigl[0,T\bigr].\) For the transformed problem, the coefficients are $R^{HK}_{1}=\widetilde{R}^{HK}_{1}=\frac{1}{2}$,  $R^{HK}_{2}=\widetilde{R}^{HK}_{2}=1$,  $Q^{HK}_{\ell}-\frac{(S^{HK}_{\ell})^2}{R^{HK}_{\ell}}=1,(\ell\in\{1,2\})$, $G^{HK}_{1}=\widetilde{G}^{HK}_{1}=\frac{1}{2}$, $G^{HK}_{2}=\widetilde{G}^{HK}_{2}=0$,  $\widetilde{Q}^{HK}_{1}-\frac{(\widetilde{S}^{HK}_{1})^2}{\widetilde{R}^{HK}_{1}}=5$, and $\widetilde{Q}^{HK}_{2}-\frac{(\widetilde{S}^{HK}_{2})^2}{\widetilde{R}^{HK}_{2}}=1$. Therefore,  the transformed problem satisfies assumption~\ref{ass: stand cod}
with $\delta=1/2$.

It follows from Theorem \ref{the: H bar K bar} that the original Hamiltonian system \eqref{equ: exa 2 FBSDE} admits a unique solution.  The correspondence established above therefore implies that
\eqref{eq:scalar-fbsde} admits a unique  solution \((X,p,\zeta,k)\in\mathcal S_{\mathbb F}^2([0,T];\mathbb R)\times\mathcal S_{\mathbb F}^2([0,T];\mathbb R)\times \mathcal L_{\mathbb F}^2([0,T];\mathbb R)\times \mathcal{M}^{2}_{\mathbb{F}}([0,T];\mathbb{R})\).

\section{Financial Applications}
\label{sec: fo alm}
\subsection{Conditional Mean-Variance Portfolio Selection with Liability}
We consider continuous-time portfolio selection with exogenous liability
payments in a regime-switching market. The market contains one risk-free
asset (a bank account) and $m$ risky assets (stocks).
Let $\bar W=(W^1,\ldots,W^{d^{\prime}})^\top$ be a $d^{\prime}$-dimensional standard Brownian
motion, where $d^{\prime}\geq m$, and let $W^{d^{\prime}+1}$ be an independent standard
Brownian motion. These Brownian motions are independent of the
Markov chain $\theta$. The asset prices satisfy

\begin{equation}\label{equ: alm bond}
	\begin{aligned}
		\mathrm dS_t^0=\ r_{(\theta_t,t)}S_t^0\,\mathrm dt,\qquad
		\mathrm dS_t^i=\ S_t^i\Bigl[
		\mu^i_{(\theta_t,t)}\,\mathrm dt+
		\sum_{j=1}^{d^{\prime}}\sigma^{ij}_{(\theta_t,t)}\,\mathrm dW_t^j
		\Bigr],\quad  S_0^0=s^0>0,\quad S_0^i=s_i>0,\quad 1\leq i\leq m.
	\end{aligned}
\end{equation}
Here $r_{(\ell,t)}$ is the risk-free interest rate,
$\mu_{(\ell,t)}=(\mu^1_{(\ell,t)},\ldots,\mu^m_{(\ell,t)})^\top$ is the
vector of expected stock return rates, and
$\sigma_{(\ell,t)}=(\sigma^{ij}_{(\ell,t)})_{m\times d^{\prime}}$ is the stock
volatility matrix in regime $\ell\in\mathbb D$.
These coefficients are bounded deterministic functions within each regime.
Set $M_{(\ell,t)}=\sigma_{(\ell,t)}\sigma_{(\ell,t)}^\top$ and assume
$M_{(\ell,t)}\geq\varepsilon \mathbf{I}_{m\times m}$ for some $\varepsilon>0$ and all
$(\ell,t)\in\mathbb D\times[0,T]$. Define the excess expected return vector
$b_{(\ell,t)}=\mu_{(\ell,t)}-r_{(\ell,t)}\mathbf1_m$.

Let $L_t$ denote the investor's cumulative net liability payments, satisfying
\begin{equation}\label{equ: dya one of lia}
	\begin{aligned}
		\mathrm dL_t= a_{(\theta_t,t)}\,\mathrm dt
		+v_{(\theta_t,t)}\,\mathrm dB_t^L,\quad L_0=\ l,
	\end{aligned}
\end{equation}
where $B^L$ is a one-dimensional standard Brownian motion. The liability
payment drift $a_{(\ell,\cdot)}$ and volatility $v_{(\ell,\cdot)}$ are
bounded deterministic scalar functions, with $v_{(\ell,t)}\geq0$.

The liability need not evolve independently of the risky asset prices.
For $1\leq j\leq d^{\prime}$, let $\varrho^j_{(\ell,t)}$ denote the instantaneous
correlation coefficient between $B^L$ and $W^j$, and write
$\varrho_{(\ell,t)}=(\varrho^1_{(\ell,t)},\ldots,
\varrho^{d^{\prime}}_{(\ell,t)})^\top$.
We assume that $\varrho_{(\ell,\cdot)}$ is a bounded deterministic
$\mathbb R^{d^{\prime}}$-valued function with $|\varrho_{(\ell,t)}|\leq1$, and specify
\(
B_t^L=\int_0^t\varrho_{(\theta_s,s)}^\top\,\mathrm d\bar W_s
+\int_0^t\sqrt{1-|\varrho_{(\theta_s,s)}|^2}\,
\mathrm dW_s^{d^{\prime}+1}.
\)
Consequently, the liability dynamics \eqref{equ: dya one of lia} can be written as
\begin{equation}\label{equ: alm liability}
	\begin{aligned}
		\mathrm dL_t={}a_{(\theta_t,t)}\,\mathrm dt
		+v_{(\theta_t,t)}\varrho_{(\theta_t,t)}^\top\,\mathrm d\bar W_t+v_{(\theta_t,t)}
		\sqrt{1-|\varrho_{(\theta_t,t)}|^2}\,\mathrm dW_t^{d^{\prime}+1},\quad 
		L_0={}l.
	\end{aligned}
\end{equation}
The decomposition in \eqref{equ: alm liability} includes the following special cases.
If $\varrho_{(\ell,t)}=0$ for all $(\ell,t)$, then $B^L=W^{d^{\prime}+1}$,
so the liability and stocks have independent Brownian drivers. Their
coefficients can still depend on the same market regime.
If $d^{\prime}>m$ and $|\varrho_{(\ell,t)}|=1$, the liability is driven entirely
by the market Brownian vector $\bar W$, but its Brownian exposure need not
be spanned by the stocks.
If $d^{\prime}=m$ and $|\varrho_{(\ell,t)}|=1$, the instantaneous Brownian exposure
of the liability can be fully hedged using the $m$ risky assets.

For later use, put
\begin{equation}\label{equ: alm liability coefficients}
	\begin{aligned}
		c_{(\ell,t)}=v_{(\ell,t)}\varrho_{(\ell,t)},\qquad
		c^j_{(\ell,t)}=v_{(\ell,t)}\varrho^j_{(\ell,t)},\quad 1\leq j\leq d^{\prime},\qquad 
		\eta_{(\ell,t)}=v_{(\ell,t)}\sqrt{1-|\varrho_{(\ell,t)}|^2}.
	\end{aligned}
\end{equation}

The investor starts with resources $w$ and an initial liability payment $l$,
leaving initial net wealth $x=w-l$. We use $X_t^u$ for the investment
account after liability payments. Trading is continuous, with no transaction costs,
and borrowing and short selling are allowed.
Let $u_t^i$ be the monetary amount invested in stock $i$ and
$u_t=(u_t^1,\ldots,u_t^m)^\top$. The amount in the bank account is
$X_t^u-\mathbf1_m^\top u_t$. Self-financing trading apart from the
prescribed liability payments gives 
\(\mathrm dX_t^u=
\bigl(X_t^u-\mathbf1_m^\top u_t\bigr)\frac{\mathrm dS_t^0}{S_t^0}
+\sum_{i=1}^m u_t^i\frac{\mathrm dS_t^i}{S_t^i}
-\mathrm dL_t.\) Applying \eqref{equ: alm bond}, \eqref{equ: alm liability}, and \eqref{equ: alm liability coefficients}, we derive 
\begin{equation}\label{equ: alm wealth}
	\begin{aligned}
		\mathrm dX_t^u={}
		\bigl[r_{(\theta_t,t)}X_t^u+
		b_{(\theta_t,t)}^\top u_t-a_{(\theta_t,t)}\bigr]\,\mathrm dt+\sum_{i=1}^{d^{\prime}}\bigl[\bigl(\sigma^{\cdot i}_{(\theta_t,t)}\bigr)^\top u_t-c^{i}_{(\theta_t,t)}\bigr]
		\mathrm d W^{i}_t-\eta_{(\theta_t,t)}\,\mathrm dW_t^{d^{\prime}+1},\quad 
		X_0^u={}x,
	\end{aligned}
\end{equation}where $\sigma^{\cdot i}_{(\ell,t)}$ denotes the $i$th column of
$\sigma_{(\ell,t)}$ for $1\leq i\leq d^{\prime}$.

The investor aims to minimize the following conditional mean-variance criterion:
\begin{equation}\label{equ: alm objective one}
	\begin{aligned}
		J\bigl(x,\ell_0;u\bigr)
		={}&\mathbb E\biggl\{
		\frac{\gamma_{\theta_T}}2\,
		\mathbb{VAR}^{\theta}_T\bigl[X_T^u\bigr]
		-\mathbb E_T^\theta\bigl[X_T^u\bigr]\biggr\},
	\end{aligned}
\end{equation}where \(\mathbb{VAR}^{\theta}_T[X_T^u]=\mathbb E_T^\theta[
(X_T^u-\mathbb E_T^\theta[X_T^u])^2]\).
Here $\gamma_\ell$ is the risk-aversion coefficient in regime $\ell$, and we
assume that $\gamma_\ell\geq\gamma_*>0$ for every $\ell\in\mathbb D$. The criterion penalizes wealth fluctuations within a
given regime history. Since $\gamma_{\theta_T}$ is $\mathcal F_T^\theta$-measurable, the tower property rewrites the criterion \eqref{equ: alm objective one} as 
\begin{equation}\label{equ: alm objective}
	\begin{aligned}
		J\bigl(x,\ell_0;u\bigr)
		={}\frac12\mathbb E\bigl\{
		\gamma_{\theta_T}\bigl(X_T^u\bigr)^2
		-\gamma_{\theta_T}
		\bigl(\mathbb E_T^\theta\bigl[X_T^u\bigr]\bigr)^2-2X_T^u\bigr\}.
	\end{aligned}
\end{equation}This portfolio selection problem is a special case of
\textbf{Problem CMF} in Section~\ref{sec: CMF of indef coe}. However, the  control weight vanishes, which violates assumption \ref{ass: stand cod}. We next verify the hypotheses of Theorem~\ref{the: H bar K bar} for this problem.

Introduce the quantities
\begin{equation}\label{equ: alm derived coefficients}
	\begin{aligned}
		\kappa_{(\ell,t)}
		=b_{(\ell,t)}^\top M_{(\ell,t)}^{-1}b_{(\ell,t)},\quad 
		\pi^L_{(\ell,t)}
		=M_{(\ell,t)}^{-1}\sigma_{(\ell,t)}c_{(\ell,t)},\quad 
		\varpi_{(\ell,t)}
		=v_{(\ell,t)}^2-
		c_{(\ell,t)}^\top\sigma_{(\ell,t)}^\top
		M_{(\ell,t)}^{-1}\sigma_{(\ell,t)}c_{(\ell,t)} .
	\end{aligned}
\end{equation}

For this model, the two Riccati equations \eqref{equ: Ricca one of Pro CMF}-\eqref{equ: Ricca two of Pro CMF} reduce
to
\begin{equation}\label{equ: alm Phi}
	\begin{aligned}
		\dot\Phi_{(\ell,t)}={}-\biggl[
		\bigl(2r_{(\ell,t)}-\kappa_{(\ell,t)}\bigr)\Phi_{(\ell,t)}
		+\sum_{j\in\mathbb D}\lambda^{\ell j}\Phi_{(j,t)}\biggr],\quad 
		\Phi_{(\ell,T)}={}\gamma_\ell,\quad\ell\in\mathbb D,
	\end{aligned}
\end{equation}and $\Pi_{(\ell,t)}\equiv0$.  Further, since $\Pi_{(\ell,t)}=0$ and $\widetilde L_{(\ell,t)}=0$, equation
\eqref{equ: Lya of Pro CMF} has the solution
$\varphi_t=\varphi_{(\theta_t,t)}$, where $\varphi_{(\ell,t)}$ satisfies
\begin{equation}\label{equ: alm varphi}
	\begin{aligned}
		\dot\varphi_{(\ell,t)}={}-\biggl[
		r_{(\ell,t)}\varphi_{(\ell,t)}
		+\sum_{j\in\mathbb D}\lambda^{\ell j}\varphi_{(j,t)}\biggr],\quad 
		\varphi_{(\ell,T)}={}-1,\quad\ell\in\mathbb D.
	\end{aligned}
\end{equation} 

Applying the  Feynman-Kac representation to \eqref{equ: alm Phi}, we have  
\(
\Phi_{(\ell,t)}
=\mathbb E[
\gamma_{\theta_T}
\mathrm{e}^{\int_t^T
	(2r_{(\theta_s,s)}-\kappa_{(\theta_s,s)})
	\mathrm ds}|\theta_t=\ell].
\) Let
\(
C_*=\max_{\ell\in\mathbb D}\operatorname*{ess\,sup}_{t\in[0,T]}
|2r_{(\ell,t)}-\kappa_{(\ell,t)}|.
\)
Then we have $\Phi_{(\ell,t)}\geq\gamma_*e^{-C_*T}>0$. Consequently,
\(
\Sigma_{(\ell,t)}=\widetilde\Sigma_{(\ell,t)}
=\Phi_{(\ell,t)}M_{(\ell,t)}
\geq\varepsilon\gamma_*e^{-C_*T}\mathbf{I}_{m\times m}.
\) Choose $\bar H=\Phi$ and $\bar K=0$. By  Lemma \ref{lem: exis of HK}, we know that the hypotheses of Theorem~\ref{the: H bar K bar} are satisfied. It follows from \eqref{equ: u CMF} that the optimal portfolio is 
\begin{equation}\label{equ: alm optimal portfolio}
	\begin{aligned}
		u_t^*={}&-M_{(\theta_t,t)}^{-1}b_{(\theta_t,t)}
		\bigl(X_t^*-\mathbb E_t^\theta\bigl[X_t^*\bigr]\bigr)
		+\pi^L_{(\theta_t,t)}-\frac{\varphi_{(\theta_t,t)}}{\Phi_{(\theta_t,t)}}
		M_{(\theta_t,t)}^{-1}b_{(\theta_t,t)}.
	\end{aligned}
\end{equation}The explicit cost formula \eqref{equ: exp cost of CMF} yields 	\begin{equation}\label{equ: alm optimal value}
	\begin{aligned}
		\mathscr V\bigl(x,\ell_0\bigr)
		= x\mathbb E\bigl[\varphi_{(\ell_0,0)}\bigr]
		+\mathbb E\int_0^T\Bigl[
		\varphi_{(\theta_t,t)}
		\bigl(b_{(\theta_t,t)}^\top\pi^L_{(\theta_t,t)}
		-a_{(\theta_t,t)}\bigr)
		-\frac{\kappa_{(\theta_t,t)}\varphi_{(\theta_t,t)}^2}
		{2\Phi_{(\theta_t,t)}}+\frac{\Phi_{(\theta_t,t)}}2
		\varpi_{(\theta_t,t)}
		\Bigr]\mathrm dt.
	\end{aligned}
\end{equation}

Summarizing the above analysis, we have the following conclusion.
\begin{Proposition}\label{pro: alm optimal portfolio}
	The conditional mean-variance portfolio selection problem defined by
	\eqref{equ: alm wealth} and \eqref{equ: alm objective} admits a unique optimal
	strategy given by \eqref{equ: alm optimal portfolio}. The optimal cost is given by \eqref{equ: alm optimal value}.
\end{Proposition}

\begin{Remark}\label{rem: alm interpretation}
	If the liability and stocks have independent Brownian drivers, set
	$\varrho_{(\ell,t)}=0$. Then $c_{(\ell,t)}=0$,
	$\eta_{(\ell,t)}=v_{(\ell,t)}$, $\pi^L_{(\ell,t)}=0$, and
	$\varpi_{(\ell,t)}=v_{(\ell,t)}^2$. All preceding conclusions remain valid.
\end{Remark}

We next characterize the conditional mean-variance curve generated by the optimal strategies as the   risk-aversion coefficient \(\gamma\) varies. Assume that the coefficient  is independent of the market regime, that is, $\gamma_\ell=\gamma>0$ for every $\ell\in\mathbb D$,
where $\gamma$ is a deterministic constant. In this case, the Feynman-Kac representations of
\eqref{equ: alm Phi} and \eqref{equ: alm varphi} give
\begin{equation}\label{equ: alm common gamma representations}
	\begin{aligned}
		\Phi_{(\ell,t)}
		={\gamma}\mathbb E\!\Bigl[
		\mathrm{e}^{\int_t^T
			(2r_{(\theta_s,s)}-\kappa_{(\theta_s,s)})\,\mathrm ds}\Big|\theta_t=\ell\Bigr],\quad 
		\varphi_{(\ell,t)}
		=-\mathbb E\!\Bigl[
		\mathrm{e}^{\int_t^T r_{(\theta_s,s)}\,\mathrm ds}\Big|\theta_t=\ell\Bigr].
	\end{aligned}
\end{equation}
Consequently, $\gamma\varphi_{(\ell,t)}/\Phi_{(\ell,t)}$ is independent
of $\gamma$, although $\Phi_{(\ell,t)}$ itself depends on $\gamma$. For $0\leq s\leq t\leq T$, define the factors \(\mathcal R_{s,t}
=\mathrm{e}^{\int_s^t r_{(\theta_{\tau},\tau)}\,\mathrm d\tau}\), \(
\mathcal H_{s,t}=\mathrm{e}^{\int_s^t(2r_{(\theta_{\tau},\tau)}-\kappa_{(\theta_{\tau},\tau)})\,\mathrm d\tau},\) which depend on the regime history. Set
\begin{equation*}\label{equ: alm conditional curve coefficients}
	\begin{aligned}
		\mathcal A_T={}&x\mathcal R_{0,T}
		+\int_0^T\mathcal R_{s,T}
		\bigl(b_{(\theta_s,s)}^\top\pi^L_{(\theta_s,s)}
		-a_{(\theta_s,s)}\bigr)\,\mathrm ds,\quad &\ 
		\mathcal B_T={}&-\int_0^T\mathcal R_{s,T}
		\kappa_{(\theta_s,s)}
		\frac{\gamma\varphi_{(\theta_s,s)}}{\Phi_{(\theta_s,s)}}\,\mathrm ds,\\
		\mathcal C_T={}&\int_0^T\mathcal H_{s,T}
		\varpi_{(\theta_s,s)}\,\mathrm ds,\quad &\ 
		\mathcal D_T={}&\int_0^T\mathcal H_{s,T}
		\kappa_{(\theta_s,s)}
		\biggl(\frac{\gamma\varphi_{(\theta_s,s)}}{\Phi_{(\theta_s,s)}}\biggr)^2
		\,\mathrm ds.
	\end{aligned}
\end{equation*}
These quantities are $\mathcal F_T^\theta$-measurable and independent of $\gamma$. 
\begin{Proposition}\label{pro: alm conditional curve}
	Suppose that $\gamma_\ell=\gamma>0$ for every $\ell\in\mathbb D$,
	where $\gamma$ is deterministic. The optimal wealth $X^*$ satisfies
	\begin{equation}\label{equ: alm conditional moment pair}
		\begin{aligned}
			\mathbb E_T^\theta\bigl[X_T^*\bigr]
			=\mathcal A_T+\frac{\mathcal B_T}{\gamma},\quad 
			\mathbb{VAR}^{\theta}_T\bigl[
			X_T^*\bigr]
			=\mathcal C_T+\frac{\mathcal D_T}{\gamma^2}.
		\end{aligned}
	\end{equation}
	Consequently, on  $\{\omega:\mathcal B_T>0\}$,
	\begin{equation}\label{equ: alm conditional mean variance relation}
		\mathbb{VAR}^{\theta}_T\bigl[
		X_T^*\bigr]
		=\mathcal C_T+\frac{\mathcal D_T}{\mathcal B_T^2}
		\bigl(\mathbb E_T^\theta\bigl[X_T^*\bigr]-\mathcal A_T\bigr)^2.
	\end{equation}
	On  $\{\omega:\mathcal B_T=0\}$, one has $\mathcal D_T=0$, and the conditional
	mean and variance are $\mathcal A_T$ and $\mathcal C_T$, respectively,
	for every $\gamma>0$.
\end{Proposition}

\begin{proof}
	Substituting \eqref{equ: alm optimal portfolio} into
	\eqref{equ: alm wealth} and taking the conditional expectation \(\mathbb E_t^\theta[\cdot]\)  yields
	\begin{equation}\label{equ: alm conditional mean ode}
		\begin{aligned}
			\frac{\mathrm d}{\mathrm dt}\mathbb E_t^\theta\bigl[X_t^*\bigr]
			={}r_{(\theta_t,t)}\mathbb E_t^\theta\bigl[X_t^*\bigr]
			+b_{(\theta_t,t)}^\top\pi^L_{(\theta_t,t)}-a_{(\theta_t,t)}-\frac1\gamma\kappa_{(\theta_t,t)}
			\frac{\gamma\varphi_{(\theta_t,t)}}{\Phi_{(\theta_t,t)}},\quad 
			\mathbb E_0^\theta\bigl[X_0^*\bigr]={}x.
		\end{aligned}
	\end{equation}
	Solving this linear equation yields the first identity in
	\eqref{equ: alm conditional moment pair}.
	
	Applying It\^o's formula to \((X_t^*-\mathbb E_t^\theta[X_t^*])^2\), taking the conditional expectation \(\mathbb{E}^{\theta}_t[\cdot]\), and using \eqref{equ: alm derived coefficients} yields
	\begin{equation}\label{equ: alm conditional variance ode}
		\begin{aligned}
			\frac{\mathrm{d}}{\mathrm{d}t}\mathbb{VAR}^{\theta}_t\bigl[
			X_t^*\bigr]
			={}\bigl(2r_{(\theta_t,t)}-\kappa_{(\theta_t,t)}\bigr)
			\mathbb{VAR}^{\theta}_t\bigl[
			X_t^*\bigr]+\varpi_{(\theta_t,t)}+\frac1{\gamma^2}\kappa_{(\theta_t,t)}
			\biggl(\frac{\gamma\varphi_{(\theta_t,t)}}{\Phi_{(\theta_t,t)}}\biggr)^2,\quad
			\mathbb{VAR}^{\theta}_0\bigl[
			X_0^*\bigr]={}0.
		\end{aligned}
	\end{equation}
	These conditional moment equations \eqref{equ: alm conditional mean ode}-\eqref{equ: alm conditional variance ode} can be read along each fixed regime path.  Solving \eqref{equ: alm conditional variance ode} gives the second identity in
	\eqref{equ: alm conditional moment pair}.
	On  $\{\omega:\mathcal B_T>0\}$, the first identity in
	\eqref{equ: alm conditional moment pair} implies
	\begin{equation}\label{equ: gamma iden}
		\begin{aligned}
			\frac1\gamma
			=\frac{\mathbb E_T^\theta[X_T^*]-\mathcal A_T}{\mathcal B_T}.
		\end{aligned}
	\end{equation}Substitution into the second identity eliminates $\gamma$ and proves
	\eqref{equ: alm conditional mean variance relation}.
	
	Finally, \eqref{equ: alm common gamma representations} gives
	$\Phi_{(\ell,t)}>0$ and $-\varphi_{(\ell,t)}>0$.
	Together with $\kappa_{(\ell,t)}\geq0$, this implies that
	$\mathcal B_T$ and $\mathcal D_T$ are both strictly positive on the event
	$\{\omega:\int_0^T\kappa_{(\theta_s,s)}\,\mathrm ds>0\}$
	and both vanish on its complement, proving the remaining assertion.
\end{proof}

\begin{Remark}\label{rem: alm conditional curve interpretation}
	Identity~\eqref{equ: alm conditional mean variance relation} describes
	the conditional moments generated by the optimal strategies of
	\eqref{equ: alm objective} as  $\gamma>0$
	varies. It is not a relation for arbitrary admissible strategies and
	does not establish the minimum conditional variance for every prescribed
	terminal conditional mean. Identity \eqref{equ: gamma iden} does not allow \(\gamma\) to depend on the future regime path: it is fixed at time zero.
	On $\{\mathcal B_T>0\}$, the conditional mean exceeds $\mathcal A_T$
	for every finite $\gamma>0$. The conditional mean and variance converge
	to $\mathcal A_T$ and $\mathcal C_T$, respectively, as $\gamma\to\infty$.
	Each optimal strategy is conditionally Pareto efficient in the following
	sense: no single admissible strategy can weakly increase its conditional mean
	and weakly decrease its conditional variance almost surely, with at least one
	strict improvement on an event of positive probability. Otherwise the
	expected objective would be smaller. 
\end{Remark}

\begin{Remark}\label{rem: alm Xie frontier reduction}
	Consider the single-regime case $\mathbb D=\{\ell\}$, so that
	$\theta_t\equiv\ell$ for all $t\in[0,T]$, and suppress the regime index
	in this remark. Then $\mathcal F_T^\theta$ is trivial, and the conditional
	mean and variance reduce to the unconditional mean and variance.
	Using the notation of
	equation~(31) in \cite{Xie Li Wang 2008}, the four coefficients reduce to
	\(
	\mathcal A_T=\mathcal D_1,
	\mathcal C_T=\mathcal D_2,
	\mathcal B_T=\mathcal D_T=e^K-1,
	\)
	where $\mathcal D_1$, $\mathcal D_2$, and $\varsigma_s$ are the
	quantities defined in that reference, with $K=\int_0^T\varsigma_s\,\mathrm ds$.
	Thus, for $K>0$,
	\(
	\frac{\mathcal D_T}{\mathcal B_T^2}
	=\frac{1}{e^K-1}
	=\frac{e^{-K}}{1-e^{-K}},
	\)
	which is the quadratic coefficient in their equation~(31).
	Consequently, \eqref{equ: alm conditional mean variance relation}
	recovers that equation on the branch $\mathbb E[X_T^*]>\mathcal A_T$;
	the endpoint with mean $\mathcal A_T$ and variance $\mathcal C_T$
	is recovered as $\gamma\to\infty$.
\end{Remark}

\subsection{Numerical examples}
\label{sec:combined_numerical_examples}

We provide a numerical example to illustrate the effectiveness of the proposed result. Suppose that there exists one stock whose   price is driven by \((W^1,W^2)\).  The  liability payments are driven by these Brownian motions  together with an additional independent Brownian motion $W^3$. The Markov chain takes values in $\{1,2\}$.  State $1$ represents a bull market and state $2$ a bear market. Time is
measured in years, with $T=2$ and initial net wealth $x=10$, measured in thousands of monetary units.  The transition intensities of \(\theta\) are \(\lambda^{12}=1\) and \(\lambda^{21}=0.7\).
The risk-aversion coefficient is the same deterministic constant in both regimes.
Table~\ref{tab:combined_baseline} lists the baseline parameters.
\begin{table}[htbp]
	\centering\small
	\caption{Baseline model parameters}
	\label{tab:combined_baseline}
	\begin{tabular}{lccccccc}
		\toprule
		Market Regime & $r_\ell$ & $\mu_\ell$ & $a_\ell$ & $v_\ell$
		& $\sigma_\ell$ & $\varrho_\ell^\top$ & $\gamma$ \\
		\midrule
		Bull \((\ell=1)\) &0.06&0.2&0.01&0.15&$(0.1,0.15)$&$(0.15,0.35)$ & \multirow{2}{*}{$8$} \\
		Bear \((\ell=2)\) &0.02&0.05&0.04&0.3&$(0.2,0.3)$&$(0.15,0.20)$ & \\
		\bottomrule
	\end{tabular}
\end{table}

\begin{figure}[!htbp]
	\centering
	\captionsetup[subfloat]{font=small,labelfont=normalfont,justification=centering}
	\subfloat[Market regime $\theta_t$\label{fig:combined_regime_path}]{%
		\makebox[0.5\linewidth][c]{\includegraphics[height=0.21\linewidth]{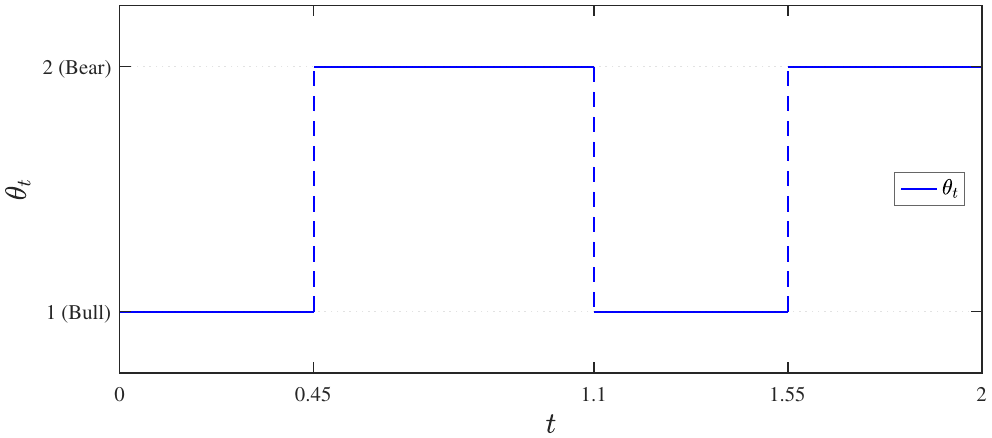}}}\hfill
	\subfloat[Optimal net wealth and conditional expectation\label{fig:combined_wealth_mean}]{%
		\makebox[0.5\linewidth][c]{\includegraphics[height=0.21\linewidth]{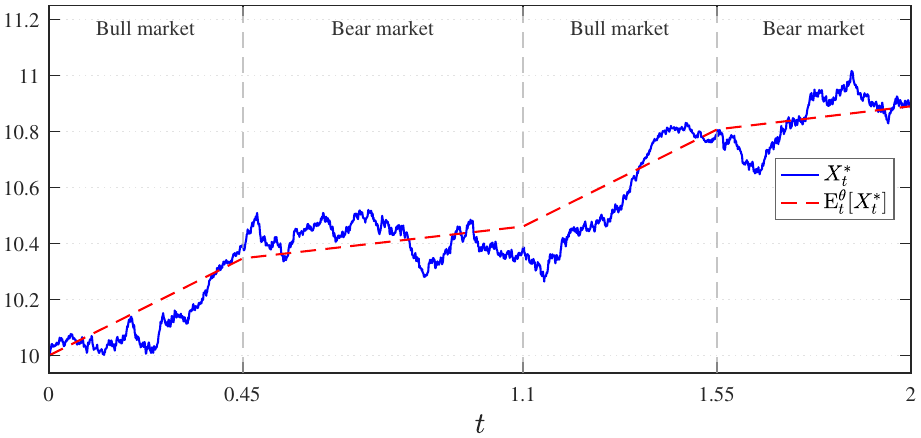}}}\\[3mm]
	\subfloat[Optimal amount invested in the stock $u_t^*$\label{fig:combined_strategy}]{%
		\makebox[0.5\linewidth][c]{\includegraphics[height=0.21\linewidth]{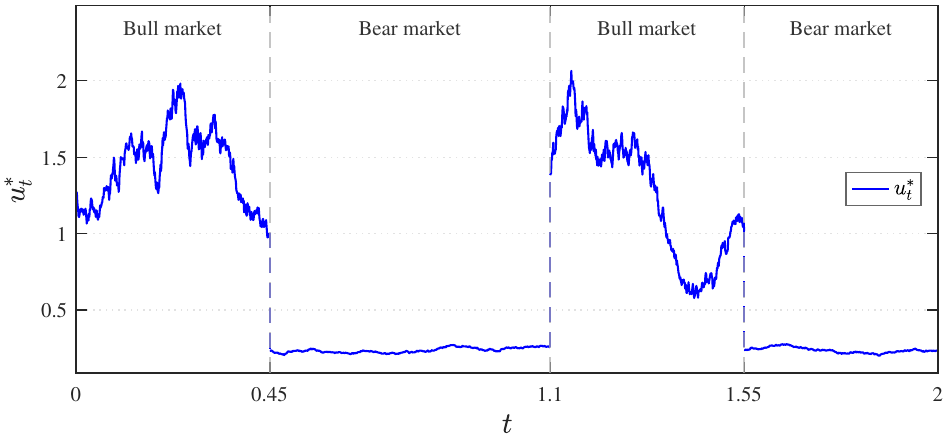}}}\hfill
	\subfloat[Optimal cost \(\mathscr V(x,\ell_0)\) by initial regime\label{fig:combined_optimal_cost}]{%
		\makebox[0.5\linewidth][c]{\includegraphics[height=0.21\linewidth]{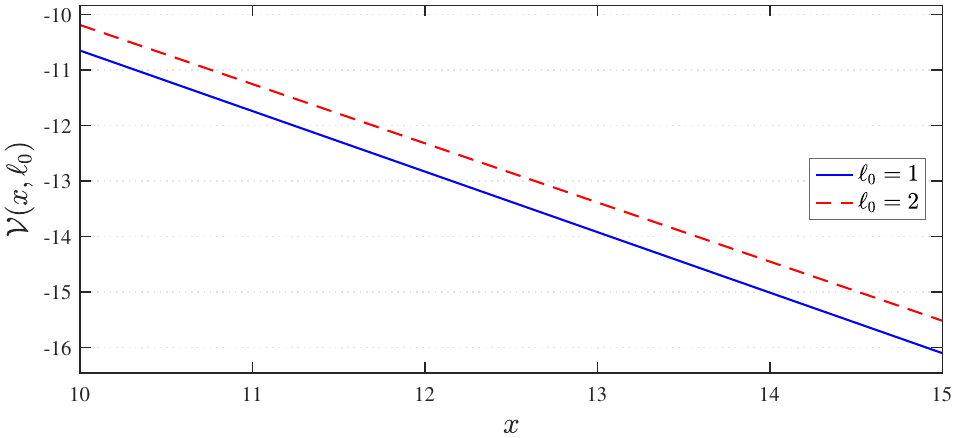}}}
	\caption{Market regime, optimal net wealth, amount invested in the stock, and corresponding cost}
	\label{fig:combined_paths_and_cost}
\end{figure}

\begin{figure}[!htb]
	\centering
	\captionsetup[subfloat]{font=small,labelfont=normalfont,justification=centering}
	\subfloat[Four cases\label{fig:combined_four_cases}]{%
		\makebox[0.33\linewidth][c]{\includegraphics[height=0.264\linewidth]{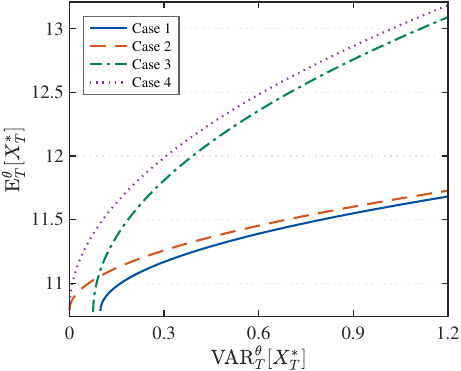}}}\hfill
	\subfloat[Initial net wealth \(x\)\label{fig:combined_initial_wealth}]{%
		\makebox[0.33\linewidth][c]{\includegraphics[height=0.264\linewidth]{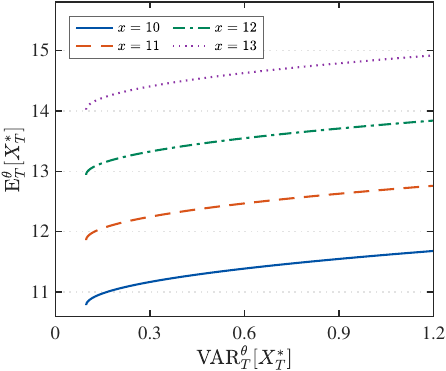}}}\hfill
	\subfloat[Liability payment drift \(a\) \label{fig:combined_liability_drift}]{%
		\makebox[0.33\linewidth][c]{\includegraphics[height=0.264\linewidth]{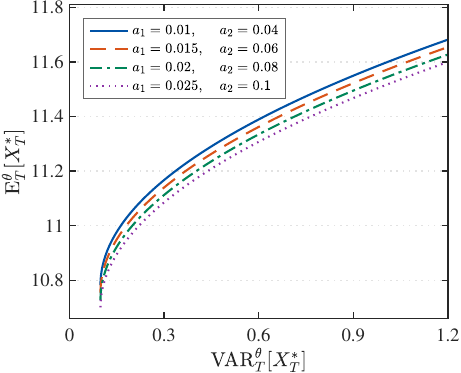}}}\\[4mm]
	\subfloat[Liability volatility \(v\)\label{fig:combined_liability_volatility}]{%
		\makebox[0.33\linewidth][c]{\includegraphics[height=0.264\linewidth]{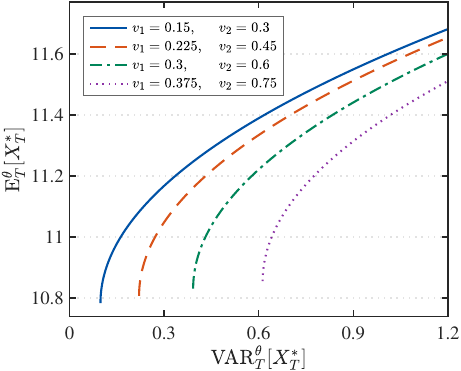}}}\hfill
	\subfloat[Correlation vector \(\varrho\)\label{fig:combined_liability_correlation}]{%
		\makebox[0.33\linewidth][c]{\includegraphics[height=0.264\linewidth]{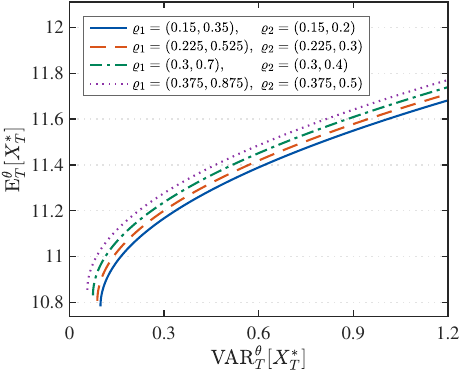}}}\hfill
	\subfloat[Expected stock return rate \(\mu\)\label{fig:combined_stock_return}]{%
		\makebox[0.33\linewidth][c]{\includegraphics[height=0.264\linewidth]{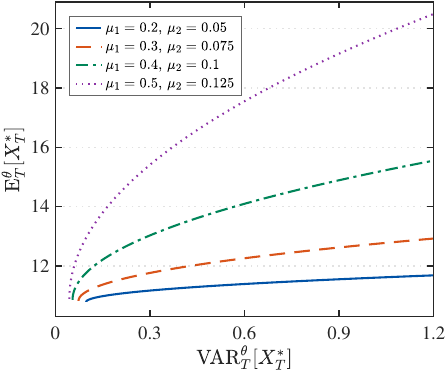}}}
	\caption{Conditional mean-variance curves for the prescribed regime history: (a) four cases in
		Table~\ref{tab:combined_cases}; (b)-(f) parameter sensitivities around
		Case 1. Blue solid curves mark the baseline.}
	\label{fig:combined_frontiers}
\end{figure}

Figure~\ref{fig:combined_regime_path} displays the path of the Markov chain \(\theta\). The trajectories in Figures~\ref{fig:combined_wealth_mean}-\ref{fig:combined_strategy} and conditional moment curves in Figure~\ref{fig:combined_frontiers} are evaluated along this   regime history. Figure~\ref{fig:combined_wealth_mean} shows a realization of optimal net wealth $X_t^*$ and its conditional expectation
$\mathbb E_t^\theta[X_t^*]$. For the displayed realization, the bull intervals show faster net growth and smaller fluctuations than
the bear intervals. Figure~\ref{fig:combined_strategy} shows the corresponding stock holding. It can be seen that  the investor holds more stock in the bull regime to pursue higher expected gains. At the switches to the bear regime, the
stock holding falls and the bank-account holding  rises.
Although $X_t^*$ and $\mathbb E_t^\theta[X_t^*]$ are continuous, the
coefficients and additive terms in \eqref{equ: alm optimal portfolio}
depend on the regime. This explains the jumps in $u_t^*$ at the
displayed switching times.  To compare the two initial regimes, we vary the
deterministic initial net wealth over $x\in[10,15]$ and compute the optimal cost $\mathscr V(x,\ell_0)$ corresponding to \(\ell_0=1\) and \(\ell_0=2\), respectively.  All the other coefficients remain as in Table~\ref{tab:combined_baseline}. Since $\varphi_{(\ell,0)}<0$ by
\eqref{equ: alm common gamma representations}, both costs decrease with
initial net wealth. As shown in Figure~\ref{fig:combined_optimal_cost}, starting
in the bull market gives a lower optimal cost than starting in the bear market.

We next compare the conditional mean-variance curves under
the optimal strategies. We consider four stock and
liability specifications, followed by a sensitivity analysis with respect to $x$, $a$, $v$,
$\varrho$, and $\mu$.  Table \ref{tab:combined_cases} summarizes the settings for the four cases, varying  the number of  Brownian motion 
factors and the presence of liability while keeping $m=1$. Unlisted market and liability coefficients retain their baseline values in Table~\ref{tab:combined_baseline}.
\begin{table}[htbp]
	\centering\small
	\caption{Coefficients of Four Different Cases}
	\label{tab:combined_cases}
	\begin{tabular}{llcl}
		\toprule
		Case & Liability & $d^{\prime}$ & Coefficients \\
		\midrule
		Case 1 & Present &2&Market and liability coefficients as in Table~\ref{tab:combined_baseline} \\
		Case 2 &Absent&2&$a_\ell=v_\ell=0$, $\varrho^{\top}_{\ell}=(0,0)$\\
		Case 3 &Present &1&$\sigma_1=0.1$, $\sigma_2=0.2$, $\varrho_{1}=\varrho_{2}=0.15$\\
		Case 4 &Absent&1&$\sigma_1=0.1$, $\sigma_2=0.2$, $a_\ell=v_\ell=\varrho_{\ell}=0$\\
		\bottomrule
	\end{tabular}
\end{table}


Figure~\ref{fig:combined_four_cases} shows the branches with $\mathbb E_T^\theta[X_T^*]>\mathcal A_T$.
Removing liability leaves $\mathcal B_T$ and $\mathcal D_T$ unchanged
within each market specification. At the common conditional terminal mean $\mathbb E_T^\theta[X_T^*]=11$,
the four curves give conditional variances approximately
$(0.1623,\ 0.0604,\ 0.0857,\ 0.0093)$.
Relative to Cases~1 and 2, Cases~3 and 4 have lower variances by
$47.1946\%$ and $84.6838\%$, respectively. These comparisons reflect the
specified changes in stock volatility and liability exposure on this
history. At this common conditional terminal mean, removing liability reduces
the conditional variance by $62.7627\%$ in Case~2 relative to Case~1 corresponding to \(d^{\prime}=2\). When \(d^{\prime}=1\),
the associated reduction in Case~4 relative to Case~3 is $89.1993\%$.

Figures~\ref{fig:combined_initial_wealth}-\ref{fig:combined_stock_return}  illustrate the effects
of initial net wealth \(x\), liability parameters \(a\), \(v\), \(\varrho\),
and the expected stock return rate \(\mu\) on the conditional mean-variance
curves. Each panel varies one parameter
group.  Figure~\ref{fig:combined_initial_wealth} shows that increasing initial net
wealth shifts the curves upwards, yielding higher conditional expected
terminal net wealth at the same conditional variance.
In Figure~\ref{fig:combined_liability_drift}, a higher liability payment drift
$a$ shifts the curves downwards because larger expected liability payments
reduce terminal net wealth. Figure~\ref{fig:combined_liability_volatility}
shows that increasing liability volatility $v$ raises the variance floor
and lowers the conditional expectation at each common admissible variance
over the displayed range, reflecting greater residual liability risk.
In Figure~\ref{fig:combined_liability_correlation},  increasing
the  correlation vectors $\varrho$ along their fixed
directions improves liability hedging, reduces residual risk and raises
the conditional expectation at a given conditional variance.
Figure~\ref{fig:combined_stock_return} shows that a higher expected
stock return rate $\mu$, with stock volatility held fixed, allows higher
conditional expected terminal net wealth at the same conditional variance.

\section{Concluding remarks}\label{sec: con rema}

This paper studies a conditional mean-field type LQ control problem for regime-switching diffusions with indefinite weights using the equivalent cost functional approach.  The quadratic adjustments for the mean and centered state in~\cite{LiLiYu2020,WangWang2024} are adapted to the conditional mean-field setting, with regime-dependent matrices accounting for transition-rate coupling and corrections to the adjoint jump integrands.  When the transformed weights in the new cost  satisfy the standard definiteness conditions, this construction yields unique solvability of the Hamiltonian system, a feedback representation of the unique optimal control, and an explicit formula for the optimal cost. The representation is determined by two coupled systems of Riccati equations and a linear BSDE driven by the Markov chain.  Two  examples are provided to illustrate our proposed results. The first example  derives the  optimal control and associated cost in a conditional mean-field type LQ problem even though the control weight is negative definite. The second one shows that our method can be used to establish the solvability of  fully coupled  CMF-FBSDEs that may not satisfy the solvability conditions imposed in the existing literature.  We also study a conditional mean-variance portfolio selection problem with exogenous liability in an incomplete regime-switching market.   The optimal investment strategy is derived, and the relation between the conditional mean and
variance of the optimal terminal net wealth is obtained.

Future work may investigate indefinite conditional mean-field type LQ control problems with random
coefficients and partial observation.




\appendix
\section{Proof of Lemma~\ref{lem: unique mini}}
\label{app: proof unique mini}

\begin{proof}We prove this conclusion in the following four steps.

	\textit{Step 1: Show that $u\rightarrow J(\xi,\ell_0; u)$ is uniformly convex on $\mathcal{U}$.}

	In this step, we use the following notation for convenience:
	\[
	\Theta_{(\ell,t)}:=R^{-1}_{(\ell,t)}S^{\top}_{(\ell,t)},
	\qquad
	\widetilde\Theta_{(\ell,t)}:=\widetilde{R}^{-1}_{(\ell,t)}\widetilde{S}^{\top}_{(\ell,t)},
	\qquad \ell\in\mathbb{D}\text{ and }t\in\bigl[0,T\bigr].
	\]

	For $u\in\mathcal U$, let $Z^u$ be the solution of the following homogeneous equation with zero initial value:
	\begin{equation}\label{equ: var X}
		\left\{\begin{aligned}
			\mathrm dZ_t^{u}={}&\bigl[A Z_t^{u}
			+\widehat A\mathbb E_t^\theta\bigl[Z_t^{u}\bigr]
			+B u_t+\widehat B\mathbb E_t^\theta\bigl[u_t\bigr]\bigr]\,\mathrm dt
			+\sum_{i=1}^N\bigl[C^iZ_t^{u}
			+\widehat C^i\mathbb E_t^\theta\bigl[Z_t^{u}\bigr]
			+D^i u_t+\widehat D^i\mathbb E_t^\theta\bigl[u_t\bigr]\bigr]
			\,\mathrm dW_t^i,\\
			Z_0^{u}={}&0.
		\end{aligned}\right.
	\end{equation}
	Define the homogeneous quadratic functional
	\begin{equation*}\label{equ: J 0}
		\begin{aligned}
			J^{0}\bigl(0,\ell_0; u\bigr)=&\frac12\mathbb E\Bigl\{
			\int_0^T\Bigl[|Z_t^{u}|_Q^2
			+|\mathbb E_t^\theta\bigl[Z_t^{u}\bigr]|_{\widehat Q}^2
			+|u_t|_R^2+|\mathbb E_t^\theta\bigl[u_t\bigr]|_{\widehat R}^2
			+2\langle S^\top Z_t^{u},u_t\rangle
			+2\langle\widehat S^\top\mathbb E_t^\theta\bigl[Z_t^{u}\bigr],
			\mathbb E_t^\theta\bigl[u_t\bigr]\rangle\Bigr]\,\mathrm dt\\
			&+|Z_T^{u}|_{G_{\theta_T}}^2
			+|\mathbb E_T^\theta\bigl[Z_T^{u}\bigr]|_{\widehat G_{\theta_T}}^2
			\Bigr\}.
		\end{aligned}
	\end{equation*}
	The mapping $u\rightarrow J(\xi,\ell_0; u)$ is uniformly convex if and only if there exists a constant $\epsilon>0$ such that
	\begin{equation}\label{equ: uni con equi}
		\begin{aligned}
			J^{0}\bigl(0,\ell_0; u\bigr)\geq \epsilon\mathbb{E}\int_{0}^{T}\left|u_{t}\right|^{2}\mathrm{d}t,\quad \forall u\in\mathcal{U}.
		\end{aligned}
	\end{equation}
	Using the properties of conditional expectation and the fact that the coefficients are adapted to $\mathbb{F}^{\theta}$, we obtain
	\begin{equation*}
		\begin{aligned}
			\mathbb E\bigl[|Z^{u}_t|_Q^2+|\mathbb{E}^{\theta}_t\bigl[Z^{u}_t\bigr]|_{\widehat Q}^2\bigr]
			&=\mathbb E\bigl[|Z^{u}_t-\mathbb{E}^{\theta}_t\bigl[Z^{u}_t\bigr]|_Q^2
			+|\mathbb{E}^{\theta}_t\bigl[Z^{u}_t\bigr]|_{\widetilde Q}^2\bigr],\\
			\mathbb E\bigl[|u_{t}|_R^2+|\mathbb{E}^{\theta}_{t}\bigl[u_{t}\bigr]|_{\widehat R}^2\bigr]
			&=\mathbb E\bigl[|u_{t}-\mathbb{E}^{\theta}_{t}\bigl[u_{t}\bigr]|_R^2
			+|\mathbb{E}^{\theta}_{t}\bigl[u_{t}\bigr]|_{\widetilde R}^2\bigr],\\
			\mathbb E\bigl[\langle S^\top Z^{u}_t,u_{t}\rangle
			+\langle\widehat S^\top\mathbb{E}^{\theta}_t\bigl[Z^{u}_t\bigr],\mathbb{E}^{\theta}_{t}\bigl[u_{t}\bigr]\rangle\bigr]
			&=\mathbb E\bigl[\langle S^\top \bigl(Z^{u}_t-\mathbb{E}^{\theta}_t\bigl[Z^{u}_t\bigr]\bigr),u_{t}-\mathbb{E}^{\theta}_{t}\bigl[u_{t}\bigr]\rangle
			+\langle\widetilde S^\top\mathbb{E}^{\theta}_t\bigl[Z^{u}_t\bigr],\mathbb{E}^{\theta}_{t}\bigl[u_{t}\bigr]\rangle\bigr].
		\end{aligned}
	\end{equation*}
	The terminal terms admit the same decomposition, with $Q,\widehat Q$ replaced by $G_{\theta_T},\widehat G_{\theta_T}$. Hence,
	\begin{equation*}
		\begin{aligned}
			J^{0}\bigl(0,\ell_0; u\bigr)=&\frac{1}{2}\mathbb E\Bigl\{
			\int_0^T\Bigl[|Z^{u}_t-\mathbb{E}^{\theta}_t\bigl[Z^{u}_t\bigr]|_Q^2
			+2\langle S^\top \bigl(Z^{u}_t-\mathbb{E}^{\theta}_t\bigl[Z^{u}_t\bigr]\bigr),u_{t}-\mathbb{E}^{\theta}_{t}\bigl[u_{t}\bigr]\rangle+|u_{t}-\mathbb{E}^{\theta}_{t}\bigl[u_{t}\bigr]|_R^2\\
			&+|\mathbb{E}^{\theta}_t\bigl[Z^{u}_t\bigr]|_{\widetilde Q}^2
			+2\langle\widetilde S^\top\mathbb{E}^{\theta}_t\bigl[Z^{u}_t\bigr],\mathbb{E}^{\theta}_{t}\bigl[u_{t}\bigr]\rangle
			+|\mathbb{E}^{\theta}_{t}\bigl[u_{t}\bigr]|_{\widetilde R}^2\Bigr]\,\mathrm dt
			+|Z^{u}_T-\mathbb{E}^{\theta}_T\bigl[Z^{u}_T\bigr]|_{G_{\theta_T}}^2
			+|\mathbb{E}^{\theta}_T\bigl[Z^{u}_T\bigr]|_{\widetilde G_{\theta_T}}^2
			\Bigr\}.
		\end{aligned}
	\end{equation*}
	For arbitrary $y\in\mathbb R^n$ and $a\in\mathbb R^m$,
	\begin{equation*}
		\begin{aligned}
			|y|_Q^2+2\langle S^\top y,a\rangle+|a|_R^2
			=|y|_{Q-S\Theta}^2+|a+\Theta y|_R^2.
		\end{aligned}
	\end{equation*}
	Applying this identity to each of the two blocks gives
	\begin{equation*}
		\begin{aligned}
			J^{0}\bigl(0,\ell_0; u\bigr)=&\frac{1}{2}\mathbb E\Bigl\{
			\int_0^T\Bigl[|Z^{u}_t-\mathbb{E}^{\theta}_t\bigl[Z^{u}_t\bigr]|_{Q-S\Theta}^2+|u_{t}-\mathbb{E}^{\theta}_{t}\bigl[u_{t}\bigr]+\Theta\bigl(Z^{u}_t-\mathbb{E}^{\theta}_t\bigl[Z^{u}_t\bigr]\bigr)|_R^2\\
			&+|\mathbb{E}^{\theta}_t\bigl[Z^{u}_t\bigr]|_{\widetilde Q-
				\widetilde S\widetilde\Theta}^2
			+|\mathbb{E}^{\theta}_{t}\bigl[u_{t}\bigr]+\widetilde\Theta\mathbb{E}^{\theta}_t\bigl[Z^{u}_t\bigr]|_{\widetilde R}^2\Bigr]\,\mathrm dt+|Z^{u}_T-\mathbb{E}^{\theta}_T\bigl[Z^{u}_T\bigr]|_{G_{\theta_T}}^2
			+|\mathbb{E}^{\theta}_T\bigl[Z^{u}_T\bigr]|_{\widetilde G_{\theta_T}}^2
			\Bigr\}.
		\end{aligned}
	\end{equation*}
	All state and terminal terms on the right-hand side are nonnegative by assumption \ref{ass: stand cod}. Therefore,
	\begin{equation}\label{equ: J 0 bound}
		\begin{aligned}
			J^{0}\bigl(0,\ell_0; u\bigr)\geq\frac{1}{2}
			\delta\,\mathbb E\int_0^T\Bigl[|u_{t}-\mathbb{E}^{\theta}_{t}\bigl[u_{t}\bigr]+\Theta\bigl(Z^{u}_t-\mathbb{E}^{\theta}_t\bigl[Z^{u}_t\bigr]\bigr)|^2+|\mathbb{E}^{\theta}_{t}\bigl[u_{t}\bigr]+\widetilde\Theta\mathbb{E}^{\theta}_t\bigl[Z^{u}_t\bigr]|^2\Bigr]\,\mathrm dt.
		\end{aligned}
	\end{equation}
	
	Define the operator $\mathscr L:\mathcal U\to\mathcal U$ by
	\begin{equation}\label{equ: ope L}
		\begin{aligned}
			\bigl(\mathscr Lu\bigr)_t
			:={}&u_t+\Theta
			\bigl(Z_t^u-\mathbb E_t^\theta\bigl[Z_t^u\bigr]\bigr)
			+\widetilde\Theta\mathbb{E}^{\theta}_t\bigl[Z^{u}_t\bigr].\\
		\end{aligned}
	\end{equation}
	Combining \eqref{equ: J 0 bound} and \eqref{equ: ope L} with the properties of conditional expectation yields
	\begin{equation}\label{equ: J0 bound}
		\begin{aligned}
			J^{0}\bigl(0,\ell_0; u\bigr)\geq\frac{1}{2}\delta\mathbb{E}\int_{0}^{T}\left|\bigl(\mathscr Lu\bigr)_t\right|^{2}\mathrm{d}t.
		\end{aligned}
	\end{equation}
	The operator $\mathscr L$ is bounded and linear. To show that $\mathscr L$ is bijective, we verify that $\mathscr L$ is both surjective and injective.
	
	Let $v\in\mathcal U$ be arbitrary, and consider the following equation for $\check{Z}^v$:
	\begin{equation}\label{equ: equ for inv ope}
		\left\{\begin{aligned}
			\mathrm d\check{Z}_t^v={}&\bigl[
			\bigl(A-B\Theta\bigr)\bigl(\check{Z}_t^v-\mathbb E_t^\theta\bigl[\check{Z}_t^v\bigr]\bigr)
			+\bigl(\widetilde A-\widetilde B\widetilde\Theta\bigr)
			\mathbb E_t^\theta\bigl[\check{Z}_t^v\bigr]+B\bigl(v_t-\mathbb E_t^\theta\bigl[v_t\bigr]\bigr)+\widetilde B\mathbb E_t^\theta\bigl[v_t\bigr]\bigr]
			\,\mathrm dt\\
			&+\sum_{i=1}^N\bigl[
			\bigl(C^i-D^i\Theta\bigr)\bigl(\check{Z}_t^v-\mathbb E_t^\theta\bigl[\check{Z}_t^v\bigr]\bigr)+\bigl(\widetilde C^i-\widetilde D^i\widetilde\Theta\bigr)
			\mathbb E_t^\theta\bigl[\check{Z}_t^v\bigr]
			+D^i\bigl(v_t-\mathbb E_t^\theta\bigl[v_t\bigr]\bigr)+\widetilde D^i\mathbb E_t^\theta\bigl[v_t\bigr]\bigr]
			\,\mathrm dW_t^i,\\
			\check{Z}_0^v={}&0.
		\end{aligned}\right.
	\end{equation}
	The boundedness of the coefficients ensures that \eqref{equ: equ for inv ope} admits a unique square-integrable solution. Define
	\[
	\bigl(\mathscr Tv\bigr)_t
	:=v_t-\Theta\bigl(\check{Z}_t^v-\mathbb E_t^\theta\bigl[\check{Z}_t^v\bigr]\bigr)
	-\widetilde\Theta\mathbb E_t^\theta\bigl[\check{Z}_t^v\bigr].
	\]
	This defines an admissible control.
	
	Set $u=\mathscr Tv$. Direct substitution into the drift coefficient gives
	\[
	\begin{aligned}
		&\ A\check{Z}_t^v+\widehat A\mathbb E_t^\theta\bigl[\check{Z}_t^v\bigr]+Bu_{t}+\widehat B\mathbb E_t^\theta\bigl[u_{t}\bigr]\\
		=&\ \bigl(A-B\Theta\bigr)\bigl(\check{Z}_t^v-\mathbb E_t^\theta\bigl[\check{Z}_t^v\bigr]\bigr)
		+\bigl(\widetilde A-\widetilde B\widetilde\Theta\bigr)\mathbb E_t^\theta\bigl[\check{Z}_t^v\bigr]
		+B\bigl(v_t-\mathbb E_t^\theta\bigl[v_t\bigr]\bigr)+\widetilde B\mathbb E_t^\theta\bigl[v_t\bigr].
	\end{aligned}
	\]
	For each $i$, the diffusion coefficient satisfies the corresponding identity
	\[
	\begin{aligned}
		&\ C^i\check{Z}_t^v+\widehat C^i\mathbb E_t^\theta\bigl[\check{Z}_t^v\bigr]+D^iu_{t}+\widehat D^i\mathbb E_t^\theta\bigl[u_{t}\bigr]\\
		=&\ \bigl(C^i-D^i\Theta\bigr)\bigl(\check{Z}_t^v-\mathbb E_t^\theta\bigl[\check{Z}_t^v\bigr]\bigr)
		+\bigl(\widetilde C^i-\widetilde D^i\widetilde\Theta\bigr)\mathbb E_t^\theta\bigl[\check{Z}_t^v\bigr]
		+D^i\bigl(v_t-\mathbb E_t^\theta\bigl[v_t\bigr]\bigr)+\widetilde D^i\mathbb E_t^\theta\bigl[v_t\bigr].
	\end{aligned}
	\]
	Thus, $\check{Z}^v$ solves \eqref{equ: var X} with control $u=\mathscr Tv$. By uniqueness, $Z^{\mathscr Tv}=\check{Z}^v$. It follows that
	\[
	\bigl(\mathscr L\bigl(\mathscr Tv\bigr)\bigr)_t
	=\bigl(\mathscr Tv\bigr)_t+\Theta\bigl(\check{Z}_t^v-\mathbb E_t^\theta\bigl[\check{Z}_t^v\bigr]\bigr)
	+\widetilde\Theta\mathbb E_t^\theta\bigl[\check{Z}_t^v\bigr]=v_t.
	\]
	Since $v$ was arbitrary, $\mathscr L$ is surjective.
	
	Suppose that $\mathscr Lu=0$. Then
	\[
	u_t=-\Theta\bigl(Z^{u}_t-\mathbb{E}^{\theta}_t\bigl[Z^{u}_t\bigr]\bigr)-\widetilde\Theta\mathbb{E}^{\theta}_t\bigl[Z^{u}_t\bigr],
	\qquad
	\mathbb{E}^{\theta}_t\bigl[u_t\bigr]=-\widetilde\Theta\mathbb{E}^{\theta}_t\bigl[Z^{u}_t\bigr].
	\]
	Substituting these identities into \eqref{equ: var X} shows that $Z^u$ satisfies \eqref{equ: equ for inv ope} for $\check{Z}^v$ with $v=0$. The identically zero process is a solution. Uniqueness therefore gives $Z^u=0$, and hence $u=0$. Thus, $\ker\mathscr L=\{0\}$.
	Combining injectivity with the right-inverse identity above gives $
	\mathscr L^{-1}=\mathscr T.$
	This establishes bijectivity.
	
	The space $\mathcal U$ is a Hilbert space, and $\mathscr L$ is a bounded linear bijection from $\mathcal U$ onto itself. The bounded inverse theorem and \eqref{equ: J0 bound} therefore give
	\begin{equation}\label{equ: J boud u}
		\begin{aligned}
			J^{0}\bigl(0,\ell_0; u\bigr)\geq\frac{\delta}{2}\mathbb{E}\int_{0}^{T}\left|\bigl(\mathscr Lu\bigr)_t\right|^{2}\mathrm{d}t\geq \frac{\delta}{2\|\mathscr L^{-1}\|_{\mathrm{op}}^{2}}\mathbb{E}\int_{0}^{T}\left|u_t\right|^{2}\mathrm{d}t.
		\end{aligned}
	\end{equation}
	Here the operator norm $\|\mathscr L^{-1}\|_{\mathrm{op}}$ is taken on the control Hilbert space. Taking $\epsilon=\frac{\delta}{2\|\mathscr L^{-1}\|_{\mathrm{op}}^{2}}$ in \eqref{equ: J boud u} yields \eqref{equ: uni con equi}, which establishes uniform convexity.
	
	\textit{Step 2: Show that $J(\xi,\ell_0; \cdot)$ is continuous.}\par\nopagebreak[4]
	
	Under assumption \ref{ass:coef in com dyna}, the Cauchy-Schwarz, Jensen, Doob, and Gronwall inequalities yield the following estimates for solutions of \eqref{equ: com sta}. There exists a constant $\mathbb{C}$ such that
	\begin{equation}\label{equ: L2 est of sta}
		\begin{aligned}
			\mathbb{E}\sup_{t\in[0,T]}\left|X^{u}_t\right|^{2}\leq \mathbb{C}\Bigl(
			1+\mathbb{E}\int_{0}^{T}\left|u_t\right|^2\mathrm{d}t\Bigr),\quad \forall u\in\mathcal{U}.
		\end{aligned}
	\end{equation}
	\begin{equation}\label{equ: conti depe}
		\begin{aligned}
			\mathbb{E}\sup_{t\in[0,T]}\left|X^{u}_t-X^{u^{\prime}}_t\right|^2\leq \mathbb{C}\mathbb{E}\int_{0}^{T}\left|u_t-u^{\prime}_t\right|^2\mathrm{d}t,\quad \forall u,u^{\prime}\in\mathcal{U}. 
		\end{aligned}
	\end{equation}
	In the remainder of this proof, the constant $\mathbb{C}$ (independent of the control $u$) may change from line to line. Combining \eqref{equ: conti depe} with assumption \ref{ass: coe in com cost}, we obtain
	\begin{equation*}
		\begin{aligned}
			|J\bigl(\xi,\ell_0;u\bigr)-J\bigl(\xi,\ell_0;u^{\prime}\bigr)|
			\leq \ \mathbb{C}\Bigl[1
			+\Bigl(\mathbb{E}\int_{0}^{T}|u_t|^2\mathrm{d}t\Bigr)^{\frac{1}{2}}+\Bigl(\mathbb{E}\int_{0}^{T}|u^{\prime}_t|^2\mathrm{d}t\Bigr)^{\frac{1}{2}}\Bigr]\Bigl(\mathbb{E}\int_{0}^{T}|u_t-u^{\prime}_t|^2\mathrm{d}t\Bigr)^{\frac{1}{2}},\quad  \forall u,u^{\prime}\in\mathcal U.
		\end{aligned}
	\end{equation*}
	This estimate establishes the continuity of $J(\xi,\ell_0; \cdot)$.
	
	\textit{Step 3: Show that the cost functional $J(\xi,\ell_0; \cdot)$ is coercive.}
	
	Let $X^0$ denote the state under zero control, with the original initial condition and inhomogeneous coefficients.
	Then $X^u=X^0+Z^u$, where $Z^{u}$ solves \eqref{equ: var X}. Expanding the cost gives
	\[
	J\bigl(\xi,\ell_0;u\bigr)
	=J\bigl(\xi,\ell_0;0\bigr)+\mathfrak l\bigl(\xi,\ell_0;u\bigr)+J^{0}\bigl(0,\ell_0;u\bigr),
	\]
	where
	\[
	\begin{aligned}
		\mathfrak l\bigl(\xi,\ell_0;u\bigr)=&\mathbb E\Bigl\{
		\int_0^T\Bigl[\langle QX_t^0,Z_t^u\rangle
		+\langle\widehat Q\mathbb E_t^\theta\bigl[X_t^0\bigr],\mathbb E_t^\theta\bigl[Z_t^u\bigr]\rangle
		+\langle S^\top X_t^0,u_t\rangle
		+\langle\widehat S^\top\mathbb E_t^\theta\bigl[X_t^0\bigr],\mathbb E_t^\theta\bigl[u_t\bigr]\rangle\\
		&+\langle q_t,Z_t^u\rangle
		+\langle\widehat q_t,\mathbb E_t^\theta\bigl[Z_t^u\bigr]\rangle
		+\langle r_t,u_t\rangle
		+\langle\widehat r_t,\mathbb E_t^\theta\bigl[u_t\bigr]\rangle\Bigr]\,\mathrm dt\\
		&+\langle G_{\theta_T}X_T^0,Z_T^u\rangle
		+\langle\widehat G_{\theta_T}\mathbb E_T^\theta\bigl[X_T^0\bigr],\mathbb E_T^\theta\bigl[Z_T^u\bigr]\rangle
		+\langle g,Z_T^u\rangle
		+\langle\widehat g,\mathbb E_T^\theta\bigl[Z_T^u\bigr]\rangle
		\Bigr\}.
	\end{aligned}
	\]
	Since $Z^u=X^u-X^0$, estimates \eqref{equ: L2 est of sta} and \eqref{equ: conti depe}, together with the Cauchy-Schwarz and Jensen inequalities, give
	\begin{equation}\label{equ: vari of cost}
		\begin{aligned}
			|\mathfrak l\bigl(\xi,\ell_0;u\bigr)|\leq \mathbb{C} \Bigl(\mathbb{E}\int_{0}^{T}|u_t|^{2}\mathrm{d}t\Bigr)^{\frac{1}{2}}.
		\end{aligned}
	\end{equation}
	Consequently, \eqref{equ: vari of cost} and \eqref{equ: J boud u} imply
	\begin{equation*}
		\begin{aligned}
			J\bigl(\xi,\ell_0;u\bigr)
			\geq J\bigl(\xi,\ell_0;0\bigr)-\mathbb{C}\Bigl(\mathbb{E}\int_{0}^{T}|u_t|^{2}\mathrm{d}t\Bigr)^{\frac{1}{2}}+\epsilon\mathbb{E}\int_{0}^{T}\left|u_{t}\right|^{2}\mathrm{d}t.
		\end{aligned}
	\end{equation*}
	It follows that
	\begin{equation*}
		\begin{aligned}
			J\bigl(\xi,\ell_0;u\bigr)\rightarrow \infty, \quad  \text{as} \quad  \mathbb{E}\int_{0}^{T}\left|u_{t}\right|^{2}\mathrm{d}t\rightarrow \infty.
		\end{aligned}
	\end{equation*}
	Therefore, $J(\xi,\ell_0; \cdot)$ is coercive.
	
	\textit{Step 4: Establish the existence and uniqueness of the  optimal control for \textbf{Problem CMF}.}

	By Steps 1-3, the cost functional $J(\xi,\ell_0;\cdot)$ is
	uniformly convex, continuous, and coercive on the Hilbert space
	$\mathcal U$. Therefore, Proposition 1.2 in Chapter II of
	\cite{Ekeland Temam 1976} implies that $J(\xi,\ell_0;\cdot)$
	attains its minimum at a unique control $u^*\in\mathcal U$.
\end{proof}

\section{Proof of Lemma~\ref{lem: wellpose of Ricc}}
\label{app: proof wellpose Ricc}

\begin{proof}Set
	\begin{equation*}\label{eq:block-coefficients}
		\begin{aligned}
			\mathcal{A}_t=&\ \mathrm{diag}\left\{A_{(1,t)}+\frac{1}{2}\lambda^{11}\mathbf{I}_{n\times n},\dots,A_{(d,t)}+\frac{1}{2}\lambda^{dd}\mathbf{I}_{n\times n}\right\}.
		\end{aligned}
	\end{equation*}Define the block-diagonal matrices \(\mathcal{B}_t,\mathcal{C}^i_t,\mathcal{D}^i_t\) \((i=1,\dots,N)\), \(\mathcal{S}_t,\mathcal{Q}_t,\mathcal{R}_t,\mathcal{G}\),  where 
	\(\mathcal B_t=\operatorname{diag}\left\{B_{(1,t)},\ldots,B_{(d,t)}\right\}\) and the other matrices are defined analogously.

	For $h=1,\ldots,d-1$, define the weighted cyclic shift matrix
	$\mathcal E_h\in\mathbb R^{nd\times nd}$ by its $n\times n$ blocks:
	\[
	(\mathcal E_h)_{\ell j}
	=\begin{cases}
		\sqrt{\lambda^{\ell j}}\,I_n,
		&j=1+((\ell-1+h)\bmod d),\\
		0,&\text{otherwise}.
	\end{cases}
	\]
	Let  $\mathcal P_t=\operatorname{diag}\left\{\Phi_{(1,t)},\ldots,\Phi_{(d,t)}\right\}$. Then
	\begin{equation*}\label{eq:coupling-identity}
		\sum_{h=1}^{d-1}\mathcal E_h\mathcal P_t\mathcal E_h^\top
		=\operatorname{diag}\left\{
		\sum_{j\ne1}\lambda^{1j}\Phi_{(j,t)},\ldots,
		\sum_{j\ne d}\lambda^{dj}\Phi_{(j,t)}
		\right\}.
	\end{equation*}
	When $d=1$, this sum is understood to be zero. 
	
	For each \(h=1,\ldots,d-1,
	\) define the additional coefficient matrices 
	\(
	\mathcal C^{N+h}_t=\mathcal E_h^\top,
	\mathcal D^{N+h}_t=0.
	\) Let $q=N+d-1$. Define
	\[
	\mathscr R_t(\mathcal{P}_t)=\mathcal R_t+
	\sum_{a=1}^{q}(\mathcal D^a_t)^\top \mathcal{P}_t\mathcal D^a_t,
	\qquad
	\mathscr L_t(\mathcal{P}_t)=\mathcal S_t+\mathcal{P}_t\mathcal B_t+
	\sum_{a=1}^{q}(\mathcal C^a_t)^\top \mathcal{P}_t\mathcal D^a_t.
	\]
	Then the solvability of \eqref{equ: Ricca one of Pro CMF} is equivalent to that of the following matrix Riccati equation:
	\begin{equation}\label{eq:lifted-riccati}
		\left\{
		\begin{aligned}
			-\dot{\mathcal P}_t
			={}&\mathcal A^\top_t\mathcal P_t+\mathcal P_t\mathcal A_t
			+\sum_{a=1}^{q}(\mathcal C^a_t)^\top
			\mathcal P_t\mathcal C^a_t+\mathcal Q_t-\mathscr L_t(\mathcal P_t)
			\mathscr R_t(\mathcal P_t)^{-1}
			\mathscr L_t(\mathcal P_t)^\top,\\
			\mathcal P_T={}&\mathcal G.
		\end{aligned}
		\right.
	\end{equation}
	Equation~\eqref{eq:lifted-riccati} is a classical Riccati equation. 
	Furthermore,
	\begin{equation*}\label{eq:lifted-positivity}
		\begin{aligned}
			\mathcal R_t&\ge\delta \mathrm{I}_{md\times md},\qquad \mathcal G\ge0,\quad 
			\mathcal Q_t-\mathcal S_t\mathcal R^{-1}_t\mathcal S^\top_t
			\ge0.
		\end{aligned}
	\end{equation*}It follows from Theorem 7.2 in Chapter 6 of  \cite{Yong Zhou 1999} that \eqref{eq:lifted-riccati} admits a unique positive semidefinite solution \(\mathcal P_t\in\mathbb S^{nd}\) such that \(	\mathscr R_t(\mathcal P_t)\) is uniformly  positive definite.  Therefore, Riccati equation \eqref{equ: Ricca one of Pro CMF} admits a unique  positive semidefinite solution \(\Phi_{(\ell,t)}\) such that 
	\[
	\Sigma_{(\ell,t)}
	=R_{(\ell,t)}+\sum_{i=1}^{N}D_{(\ell,t)}^{i\top}\Phi_{(\ell,t)} D_{(\ell,t)}^i
	\ge R_{(\ell,t)}\ge\delta \mathrm{I}_{m\times m}.
	\]
	
	Next, we establish the solvability of the second Riccati system.  For each $\ell\in\mathbb D$ and $t\in[0,T]$, define
	\begin{equation*}\label{eq:pi-effective-weights}
		\left\{\begin{aligned}
			\mathsf R_{(\ell,t)}
			&=\widetilde R_{(\ell,t)}+
			\sum_{i=1}^{N}\widetilde D_{(\ell,t)}^{i\top}
			\Phi_{(\ell,t)}\widetilde D_{(\ell,t)}^i
			=\widetilde\Sigma_{(\ell,t)},\\
			\mathsf S_{(\ell,t)}
			&=\widetilde S_{(\ell,t)}+
			\sum_{i=1}^{N}\widetilde C_{(\ell,t)}^{i\top}
			\Phi_{(\ell,t)}\widetilde D_{(\ell,t)}^i,\\
			\mathsf Q_{(\ell,t)}
			&=\widetilde Q_{(\ell,t)}+
			\sum_{i=1}^{N}\widetilde C_{(\ell,t)}^{i\top}
			\Phi_{(\ell,t)}\widetilde C_{(\ell,t)}^i.
		\end{aligned}\right.
	\end{equation*}
	These matrices do not depend on the unknown $\Pi_{(\ell,t)}$, and
	$\mathsf R_{(\ell,t)}\ge\widetilde R_{(\ell,t)}\ge\delta\mathbf{I}_{m\times m}$.
	Since $\widetilde L_{(\ell,t)}=\Pi_{(\ell,t)}\widetilde B_{(\ell,t)}+\mathsf S_{(\ell,t)}$,
	expanding the quadratic term in~\eqref{equ: Ricca two of Pro CMF} gives
	\begin{equation}\label{eq:pi-reduced}
		\left\{
		\begin{aligned}
			-\dot\Pi_{(\ell,t)}={}&
			\overline A_{(\ell,t)}^\top\Pi_{(\ell,t)}+\Pi_{(\ell,t)}\overline A_{(\ell,t)}
			+\overline Q_{(\ell,t)}-\Pi_{(\ell,t)}\widetilde B_{(\ell,t)}\mathsf R_{(\ell,t)}^{-1}
			\widetilde B_{(\ell,t)}^\top\Pi_{(\ell,t)}
			+\sum_{j=1}^{d}\lambda^{\ell j}\Pi_{(j,t)},\\
			\Pi_{(\ell,T)}={}&\widetilde G_\ell,
		\end{aligned}
		\right.
	\end{equation}
	where
	\begin{equation*}\label{eq:pi-reduced-coefficients}
		\begin{aligned}
			\overline A_{(\ell,t)}
			=\widetilde A_{(\ell,t)}
			-\widetilde B_{(\ell,t)}\mathsf R_{(\ell,t)}^{-1}\mathsf S_{(\ell,t)}^\top,\quad 
			\overline Q_{(\ell,t)}
			=\mathsf Q_{(\ell,t)}-\mathsf S_{(\ell,t)}\mathsf R_{(\ell,t)}^{-1}
			\mathsf S_{(\ell,t)}^\top.
		\end{aligned}
	\end{equation*}
	
	We next verify $\overline Q_{(\ell,t)}\ge0$. Set
	\[
	\mathbf P_{(\ell,t)}=\mathbf{I}_{N\times N}\otimes\Phi_{(\ell,t)},\qquad
	\mathbf C_{(\ell,t)}=\begin{pmatrix}
		\widetilde C_{(\ell,t)}^1\\ \vdots\\ \widetilde C_{(\ell,t)}^N
	\end{pmatrix},\qquad
	\mathbf D_{(\ell,t)}=\begin{pmatrix}
		\widetilde D_{(\ell,t)}^1\\ \vdots\\ \widetilde D_{(\ell,t)}^N
	\end{pmatrix}.
	\]
	Thus $\mathbf P_{(\ell,t)}\ge0$ and
	$\mathsf R_{(\ell,t)}=\widetilde R_{(\ell,t)}+
	\mathbf D_{(\ell,t)}^\top\mathbf P_{(\ell,t)}\mathbf D_{(\ell,t)}$.
	Define
	\[
	\begin{aligned}
		F_{(\ell,t)}=\widetilde R_{(\ell,t)}^{-1}\widetilde S_{(\ell,t)}^\top,\quad 
		\mathbf C_{(\ell,t)}^0=\mathbf C_{(\ell,t)}-\mathbf D_{(\ell,t)} F_{(\ell,t)},\quad 
		Q_{(\ell,t)}^0=\widetilde Q_{(\ell,t)}
		-\widetilde S_{(\ell,t)}\widetilde R_{(\ell,t)}^{-1}
		\widetilde S_{(\ell,t)}^\top\ge0.
	\end{aligned}
	\]
	Writing $H_{(\ell,t)}=(\mathbf C_{(\ell,t)}^0)^\top
	\mathbf P_{(\ell,t)}\mathbf D_{(\ell,t)}$, we obtain 
	\[
	\left\{\begin{aligned}
		\mathsf S_{(\ell,t)}&=H_{(\ell,t)}+F_{(\ell,t)}^\top\mathsf R_{(\ell,t)},\\
		\mathsf Q_{(\ell,t)}&=Q_{(\ell,t)}^0
		+(\mathbf C_{(\ell,t)}^0)^\top\mathbf P_{(\ell,t)}\mathbf C_{(\ell,t)}^0+H_{(\ell,t)} F_{(\ell,t)}+F_{(\ell,t)}^\top H_{(\ell,t)}^\top
		+F_{(\ell,t)}^\top\mathsf R_{(\ell,t)} F_{(\ell,t)}.
	\end{aligned}\right.
	\]
	Consequently,
	\begin{equation}\label{eq:pi-weight-decomposition}
		\begin{aligned}
			\overline Q_{(\ell,t)}
			={}&Q_{(\ell,t)}^0+(\mathbf C_{(\ell,t)}^0)^\top
			\Bigl[\mathbf P_{(\ell,t)}-\mathbf P_{(\ell,t)}\mathbf D_{(\ell,t)}\mathsf R_{(\ell,t)}^{-1}
			\mathbf D_{(\ell,t)}^\top\mathbf P_{(\ell,t)}\Bigr]
			\mathbf C_{(\ell,t)}^0.
		\end{aligned}
	\end{equation}
	Let $\Gamma_{(\ell,t)}=\mathbf P_{(\ell,t)}^{1/2}
	\mathbf D_{(\ell,t)}\widetilde R_{(\ell,t)}^{-1/2}$. Since
	\[
	\mathsf R_{(\ell,t)}=\widetilde R_{(\ell,t)}^{1/2}
	(\mathbf{I}_{m\times m}+\Gamma_{(\ell,t)}^\top\Gamma_{(\ell,t)})
	\widetilde R_{(\ell,t)}^{1/2},
	\]
	the identity
	\[
	\mathbf{I}_{Nn\times Nn}-\Gamma_{(\ell,t)}
	(\mathbf{I}_{m\times m}+\Gamma_{(\ell,t)}^\top\Gamma_{(\ell,t)})^{-1}\Gamma_{(\ell,t)}^\top
	=(\mathbf{I}_{Nn\times Nn}+\Gamma_{(\ell,t)}\Gamma_{(\ell,t)}^\top)^{-1}
	\]yields
	\begin{equation}\label{eq:pi-square-root-identity}
		\begin{aligned}
			\mathbf P_{(\ell,t)}-\mathbf P_{(\ell,t)}\mathbf D_{(\ell,t)}
			\mathsf R_{(\ell,t)}^{-1}\mathbf D_{(\ell,t)}^\top\mathbf P_{(\ell,t)}
			=&\ \mathbf P_{(\ell,t)}^{1/2}
			\left[\mathbf{I}_{Nn\times Nn}-\Gamma_{(\ell,t)}
			(\mathbf{I}_{m\times m}+\Gamma_{(\ell,t)}^\top\Gamma_{(\ell,t)})^{-1}\Gamma_{(\ell,t)}^\top
			\right]\mathbf P_{(\ell,t)}^{1/2}\\
			=&\ \mathbf P_{(\ell,t)}^{1/2}
			(\mathbf{I}_{Nn\times Nn}+\Gamma_{(\ell,t)}\Gamma_{(\ell,t)}^\top)^{-1}
			\mathbf P_{(\ell,t)}^{1/2}\ge0.
		\end{aligned}
	\end{equation}
	Combining~\eqref{eq:pi-weight-decomposition} with
	\eqref{eq:pi-square-root-identity} and $Q_{(\ell,t)}^0\ge0$, we obtain
	\begin{equation*}\label{eq:pi-positivity}
		\overline Q_{(\ell,t)}\ge0,\qquad
		\mathsf R_{(\ell,t)}\ge\delta\mathbf{I}_{m\times m},\qquad
		\widetilde G_\ell\ge0.
	\end{equation*}Moreover, all coefficients in~\eqref{eq:pi-reduced} are bounded, since $\Phi_{(\ell,t)}$ is bounded and
	$\|\mathsf R_{(\ell,t)}^{-1}\|\le\delta^{-1}$.

	Equation~\eqref{eq:pi-reduced} is a standard coupled Riccati system.
	Thus the solvability argument for the first Riccati equation applies. This proves that   \eqref{equ: Ricca two of Pro CMF} admits a unique positive semidefinite solution \(\Pi_{(\ell,t)}\) for any \(\ell\in\mathbb{D}\).

	Once the Riccati equations are solved, the driver in  \eqref{equ: Lya of Pro CMF} has a bounded linear coefficient and a square-integrable inhomogeneous term. The terminal value is also square-integrable. The  existence and uniqueness of solutions to  \eqref{equ: Lya of Pro CMF} follow from Theorem~6.2 in \cite{Cohen  Elliott 2008}.
\end{proof}

\section{Proof of Theorem~\ref{the: main result definite}}
\label{app: proof main definite}

\begin{proof}
	Let $\Phi$, $\Pi$, and $(\varphi,\rho)$ be the solutions provided by Lemma~\ref{lem: wellpose of Ricc}. It can be easily verified that  equation \eqref{equ: clos-loo sys CMF} admits a unique solution $X^*\in\mathcal S_{\mathbb F}^2([0,T];\mathbb R^n)$, and \eqref{equ: u CMF} defines an admissible control $u^*\in\mathcal U$.  Taking conditional expectation  \(\mathbb{E}^{\theta}_t[\cdot]\) on both sides of \eqref{equ: clos-loo sys CMF} and applying Lemma 2.4 in \cite{Lv Xiong Zhang 2023} yields 
	\begin{equation*}
		\left\{\begin{aligned}
			\mathrm{d}\mathbb{E}^{\theta}_t[X^{*}_t]=&\ \bigl[\bigl(\widetilde{A}-\widetilde{B} \widetilde{\Sigma}^{-1}\widetilde{L}^{\top}\bigr)\mathbb{E}^{\theta}_t\bigl[X^{*}_t\bigr]+\alpha-\widetilde{B}\widetilde{\Sigma}^{-1}\hbar\bigr]\mathrm{d}t,\\
			\mathbb{E}^{\theta}_0[X^{*}_0]=&\ \mathbb{E}\xi.
		\end{aligned}\right.
	\end{equation*}Applying It\^o's formula to the expression for \(p\) in \eqref{equ: p q k CMF} and using \eqref{equ: Ricca one of Pro CMF}-\eqref{equ: Lya of Pro CMF}, we verify the adjoint equation and its terminal condition.  Direct substitution verifies the state equation and the stationarity condition in \eqref{equ: Ham for com}. Thus, $(X^*,u^*,p,\zeta,k)$ solves \eqref{equ: Ham for com}, and Lemma~\ref{lem: nece and suff of optim} establishes the optimality of $(X^*,u^*)$.
	
	It remains to establish the optimal cost formula \eqref{equ: exp cost of CMF}. For an arbitrary \(u\in\mathcal U\), let \(X^u\) be the corresponding state process. Applying It\^o's formula to
	\(	\langle \bigl(X^{u}_t-\mathbb{E}^{\theta}_t\bigl[X^{u}_t\bigr]\bigr), \Phi_{(\theta_t,t)}\bigl(X^{u}_t-\mathbb{E}^{\theta}_t\bigl[X^{u}_t\bigr]\bigr)\rangle+\langle \mathbb{E}^{\theta}_t\bigl[X^{u}_t\bigr], \Pi_{(\theta_t,t)}\mathbb{E}^{\theta}_t\bigl[X^{u}_t\bigr]\rangle+2\langle\mathbb{E}^{\theta}_t\bigl[X^{u}_t\bigr],\varphi_{t}\rangle,\) integrating from $0$ to $T$, and taking expectations yield
	\begin{equation}\label{equ: ito for exp cos}
		\begin{aligned}
			&\ \mathbb{E}\bigl[|X^{u}_{T}|^2_{G_{\theta_T}}+|\mathbb E_T^\theta\bigl[X^{u}_{T}\bigr]|^2_{\widehat{G}_{\theta_T}}+2\langle X^{u}_{T},g\rangle+2\langle\mathbb E_T^\theta\bigl[X^{u}_{T}\bigr],\widehat{g}\rangle\bigr]\\
			=&\ \mathbb{E}\Bigl\{\int_{0}^{T}\Bigl[|u_t-\mathbb E_t^\theta\bigl[u_t\bigr]|^{2}_{\sum_{i=1}^{N}D^{i\top}\Phi D^{i}}+2\langle\bigl(u_t-\mathbb E_t^\theta\bigl[u_t\bigr]\bigr),\bigl(L-S\bigr)^{\top}\bigl(X^{u}_{t}-\mathbb E_t^\theta\bigl[X^{u}_{t}\bigr]\bigr)\rangle+|\mathbb E_t^\theta\bigl[u_t\bigr]|^{2}_{\sum_{i=1}^{N}\widetilde{D}^{i\top}\Phi \widetilde{D}^{i}}\\
			&\ +2\langle\mathbb E_t^\theta\bigl[u_t\bigr],\bigl(\widetilde{L}-\widetilde{S}\bigr)^{\top}\mathbb E_t^\theta\bigl[X^{u}_{t}\bigr]+\hbar-\widetilde{r}\rangle+|X^{u}_{t}-\mathbb E_t^\theta\bigl[X^{u}_{t}\bigr]|^2_{L\Sigma^{-1}L^{\top}-Q}+|\mathbb E_t^\theta\bigl[X^{u}_{t}\bigr]|^2_{\widetilde{L}\widetilde{\Sigma}^{-1}\widetilde{L}^{\top}-\widetilde{Q}}\\
			&\ +2\langle\mathbb E_t^\theta\bigl[X^{u}_{t}\bigr],\widetilde{L}\widetilde{\Sigma}^{-1}\hbar-\widetilde{q}\rangle+2\langle\alpha,\varphi\rangle+\sum_{i=1}^{N}\left|\beta^{i}_{t}\right|^{2}_{\Phi_{(\theta_t,t)}}\Bigr]\mathrm{d}t+|\xi-\mathbb{E}\xi|^2_{\Phi_{(\ell_0,0)}}+\left|\mathbb{E}\xi\right|^{2}_{\Pi_{(\ell_0,0)}}+2\langle\mathbb{E}\xi,\varphi_{0} \rangle\Bigr\}.
		\end{aligned}
	\end{equation}
	Substituting the right-hand side of \eqref{equ: ito for exp cos} for the terminal cost in \eqref{equ: com cos}, completing the squares in $u_t-\mathbb E_t^\theta[u_t]$ and $\mathbb E_t^\theta[u_t]$, and using the positive definiteness of $\Sigma$ and $\widetilde\Sigma$ established in Lemma~\ref{lem: wellpose of Ricc}, we obtain
	\begin{equation*}
		\begin{aligned}
			J\bigl(\xi,\ell_0; u\bigr)
			=&\ \frac{1}{2}\mathbb{E}\Bigl\{\int_{0}^{T}\Bigl[\sum_{i=1}^{N}\left|\beta^{i}_{t}\right|^{2}_{\Phi_{(\theta_t,t)}}-\left|\hbar_{t}\right|^{2}_{\widetilde{\Sigma}^{-1}_{(\theta_t,t)}}+2\alpha^{\top}_{t}\varphi_{t}\Bigr]\mathrm{d}t +\left|\xi-\mathbb{E}\xi\right|^{2}_{\Phi_{(\ell_0,0)}}+\left|\mathbb{E}\xi\right|^{2}_{\Pi_{(\ell_0,0)}}+2\langle\mathbb{E}\xi,\varphi_{0} \rangle\Bigr\}\\
			&\ +\frac12\mathbb E\Bigl\{\int_0^T\Bigl[
			\left|u_t-\mathbb E_t^\theta\bigl[u_t\bigr]+\Sigma^{-1}L^\top \bigl(X^{u}_{t}-\mathbb E_t^\theta\bigl[X^{u}_{t}\bigr]\bigr)\right|_\Sigma^2+
			\left|\mathbb E_t^\theta\bigl[u_t\bigr]+\widetilde\Sigma^{-1}
			\bigl(\widetilde L^\top \mathbb E_t^\theta\bigl[X^{u}_{t}\bigr]+\hbar_t\bigr)\right|_{\widetilde\Sigma}^2\Bigr]
			\,\mathrm dt\Bigr\}\\
			\geq&\  \frac{1}{2}\mathbb{E}\Bigl\{\int_{0}^{T}\Bigl[\sum_{i=1}^{N}\left|\beta^{i}_{t}\right|^{2}_{\Phi_{(\theta_t,t)}}-\left|\hbar_{t}\right|^{2}_{\widetilde{\Sigma}^{-1}_{(\theta_t,t)}}+2\alpha^{\top}_{t}\varphi_{t}\Bigr]\mathrm{d}t +\left|\xi-\mathbb{E}\xi\right|^{2}_{\Phi_{(\ell_0,0)}}+\left|\mathbb{E}\xi\right|^{2}_{\Pi_{(\ell_0,0)}}+2\langle\mathbb{E}\xi,\varphi_{0} \rangle\Bigr\}.\\
		\end{aligned}
	\end{equation*}For \(u=u^*\), \eqref{equ: u CMF} and its conditional expectation imply that both square terms vanish. Hence, the optimal cost is given by \eqref{equ: exp cost of CMF}.
\end{proof}

\section{Proof of Theorem~\ref{the: H bar K bar}}
\label{app: proof H bar K bar}

\begin{proof}
	Fix $\bar{H},\bar{K}\in\Delta$ such that
	\textbf{Problem CMF$^{\bar{H}\bar{K}}$}
	satisfies assumption~\ref{ass: stand cod}.
	By Lemma~\ref{lem: well-pose of FBSDE}, the transformed Hamiltonian system \eqref{equ: Ham for HK} admits a unique solution
	$(X^{\bar{H}\bar{K}},u^{\bar{H}\bar{K}},p^{\bar{H}\bar{K}},\zeta^{\bar{H}\bar{K}},k^{\bar{H}\bar{K}})$.
	Lemma~\ref{lem: Ham of indefi} then implies that the original Hamiltonian system \eqref{equ: Ham for com} also admits a unique solution.
	By Lemma~\ref{lem: wellpose of Ricc}, the transformed Riccati equations \eqref{equ: Ricca one of Pro CMF HK} and \eqref{equ: Ricca two of Pro CMF HK} admit
	unique positive semidefinite solutions
	$\Phi^{\bar{H}\bar{K}}$ and
	$\Pi^{\bar{H}\bar{K}}$.
	The associated BSDE \eqref{equ: Lya of Pro CMF HK} admits a unique solution
	$(\varphi^{\bar{H}\bar{K}},
	\rho^{\bar{H}\bar{K}})$.
	
	Define
	\begin{equation}\label{equ: Phi Pi HKbar}
		\begin{aligned}
			\Phi_{(\ell,t)}
			=\Phi^{\bar{H}\bar{K}}_{(\ell,t)}
			+\bar{H}_{(\ell,t)},
			\qquad
			\Pi_{(\ell,t)}
			=\Pi^{\bar{H}\bar{K}}_{(\ell,t)}
			+\bar{K}_{(\ell,t)},\quad \text{and} \quad \bigl(\varphi,\rho\bigr)
			=\bigl(\varphi^{\bar{H}\bar{K}},
			\rho^{\bar{H}\bar{K}}\bigr).
		\end{aligned}
	\end{equation}
	By Lemma~\ref{lem: Ricca of indefi}, these functions and processes
	solve \eqref{equ: Ricca one of Pro CMF}-\eqref{equ: Lya of Pro CMF}.
	
	Using \eqref{equ: coe of J HK}, we obtain
	\begin{equation}\label{equ:equa of Sigma bar HK}
		\left\{\begin{aligned}
			\Sigma_{(\ell,t)}
			&=
			\Sigma^{\bar{H}\bar{K}}_{(\ell,t)}
			=
			R^{\bar{H}\bar{K}}_{(\ell,t)}
			+\sum_{i=1}^{N}
			D^{i\top}_{(\ell,t)}
			\Phi^{\bar{H}\bar{K}}_{(\ell,t)}
			D^i_{(\ell,t)}
			\geq \delta\mathbf{I}_{m\times m},\\
			\widetilde{\Sigma}_{(\ell,t)}
			&=
			\widetilde{\Sigma}^{\bar{H}\bar{K}}_{(\ell,t)}
			=
			\widetilde{R}^{\bar{H}\bar{K}}_{(\ell,t)}
			+\sum_{i=1}^{N}
			\widetilde{D}^{i\top}_{(\ell,t)}
			\Phi^{\bar{H}\bar{K}}_{(\ell,t)}
			\widetilde{D}^i_{(\ell,t)}
			\geq \delta\mathbf{I}_{m\times m}.
		\end{aligned}\right.
	\end{equation}
	These inequalities follow from
	$\Phi^{\bar{H}\bar{K}}_{(\ell,t)}\geq0$
	and assumption~\ref{ass: stand cod}
	for the transformed problem.
	
	Define 
	\(	\hbar_t^{\bar{H}\bar{K}}
	=
	\widetilde{B}_{(\theta_t,t)}^\top\varphi^{\bar{H}\bar{K}}_t
	+\widetilde{r}^{\bar{H}\bar{K}}_t
	+\sum_{i=1}^{N}
	\widetilde{D}^{i\top}_{(\theta_t,t)}
	\Phi^{\bar{H}\bar{K}}_{(\theta_t,t)}
	\beta_t^i.\)
	Combining the coefficient identities in \eqref{equ: coe of J HK} with \eqref{equ: Phi Pi HKbar} gives
	\begin{equation}\label{equ: L HKbar}
		\begin{aligned}
			L^{\bar{H}\bar{K}}_{(\ell,t)}
			=L_{(\ell,t)},
			\qquad
			\widetilde{L}^{\bar{H}\bar{K}}_{(\ell,t)}
			=\widetilde{L}_{(\ell,t)},
			\qquad
			\hbar_t^{\bar{H}\bar{K}}=\hbar_t.
		\end{aligned}
	\end{equation}
	
	Consider the following closed-loop system:
	\begin{equation}\label{equ: clos-loo sys CMF bar HK}
		\left\{\begin{aligned}
			\mathrm{d}X^{\bar{H}\bar{K}*}_t=&\ \bigl[\bigl(A-B \bigl(\Sigma^{\bar{H}\bar{K}}\bigr)^{-1}\bigl(L^{\bar{H}\bar{K}}\bigr)^{\top}\bigr)\bigl(X^{\bar{H}\bar{K}*}_t-\mathbb{E}^{\theta}_t\bigl[X^{\bar{H}\bar{K}*}_t\bigr]\bigr)\\
			&\  +\bigl(\widetilde{A}-\widetilde{B} \bigl(\widetilde{\Sigma}^{\bar{H}\bar{K}}\bigr)^{-1}\bigl(\widetilde{L}^{\bar{H}\bar{K}}\bigr)^{\top}\bigr)\mathbb{E}^{\theta}_t\bigl[X^{\bar{H}\bar{K}*}_t\bigr]+\alpha-\widetilde{B}\bigl(\widetilde{\Sigma}^{\bar{H}\bar{K}}\bigr)^{-1}\hbar^{\bar{H}\bar{K}}\bigr]\mathrm{d}t\\
			&\  +\sum_{i=1}^{N}\bigl[\bigl(C^{i}-D^{i} \bigl(\Sigma^{\bar{H}\bar{K}}\bigr)^{-1}\bigl(L^{\bar{H}\bar{K}}\bigr)^{\top}\bigr)\bigl(X^{\bar{H}\bar{K}*}_t-\mathbb{E}^{\theta}_t\bigl[X^{\bar{H}\bar{K}*}_t\bigr]\bigr)\\
			&\  +\bigl(\widetilde{C}^{i}-\widetilde{D}^{i} \bigl(\widetilde{\Sigma}^{\bar{H}\bar{K}}\bigr)^{-1}\bigl(\widetilde{L}^{\bar{H}\bar{K}}\bigr)^{\top}\bigr)\mathbb{E}^{\theta}_t\bigl[X^{\bar{H}\bar{K}*}_t\bigr]+\beta^{i}-\widetilde{D}^{i}\bigl(\widetilde{\Sigma}^{\bar{H}\bar{K}}\bigr)^{-1}\hbar^{\bar{H}\bar{K}}\bigr]\mathrm{d}W^i_t,\\
			X^{\bar{H}\bar{K}*}_0=&\ \xi,
		\end{aligned}\right.
	\end{equation}and define 
	\begin{equation}\label{equ: u CMF bar HK}
		\begin{aligned}
			u^{\bar{H}\bar{K}*}_t=&\ -\bigl(\Sigma^{\bar{H}\bar{K}}\bigr)^{-1}\bigl(L^{\bar{H}\bar{K}}\bigr)^{\top}\bigl(X^{\bar{H}\bar{K}*}_t-\mathbb{E}^{\theta}_t\bigl[X^{\bar{H}\bar{K}*}_t\bigr]\bigr)-\bigl(\widetilde{\Sigma}^{\bar{H}\bar{K}}\bigr)^{-1}\bigl(\bigl(\widetilde{L}^{\bar{H}\bar{K}}\bigr)^{\top}\mathbb{E}^{\theta}_t\bigl[X^{\bar{H}\bar{K}*}_t\bigr]+\hbar_{t}\bigr).\\
		\end{aligned}
	\end{equation}
	
	By Theorem~\ref{the: main result definite}
	and Lemma~\ref{lem: unique mini},
	the transformed problem admits a unique optimal pair
	$(X^{\bar{H}\bar{K}*},u^{\bar{H}\bar{K}*})$, which satisfies  \eqref{equ: clos-loo sys CMF bar HK}-\eqref{equ: u CMF bar HK}. 
	The identities  \eqref{equ:equa of Sigma bar HK} and \eqref{equ: L HKbar} show that the transformed closed-loop system \eqref{equ: clos-loo sys CMF bar HK} and feedback control \eqref{equ: u CMF bar HK} coincide with the state equation \eqref{equ: clos-loo sys CMF} and the feedback formula for \(u^*\) in \eqref{equ: u CMF}, respectively. Hence,
	\(
	\bigl(X^*,u^*\bigr)
	=\bigl(X^{\bar{H}\bar{K}*},u^{\bar{H}\bar{K}*}\bigr).
	\)
	
	The adjoint processes associated with this optimal pair
	are given by
	\begin{equation*}
		\left\{\begin{aligned}
			p^{\bar{H}\bar{K}}_t=&\ \Phi^{\bar{H}\bar{K}}\bigl(X^{\bar{H}\bar{K}*}_t-\mathbb{E}^{\theta}_t\bigl[X^{\bar{H}\bar{K}*}_t\bigr]\bigr)+\Pi^{\bar{H}\bar{K}}\mathbb{E}^{\theta}_t\bigl[X^{\bar{H}\bar{K}*}_t\bigr]+\varphi^{\bar{H}\bar{K}}_{t},\\
			\zeta^{\bar{H}\bar{K},i}_t=&\ \Phi^{\bar{H}\bar{K}}\bigl[\bigl(C^{i}-D^{i} \bigl(\Sigma^{\bar{H}\bar{K}}\bigr)^{-1}\bigl(L^{\bar{H}\bar{K}}\bigr)^{\top}\bigr)\bigl(X^{\bar{H}\bar{K}*}_t-\mathbb{E}^{\theta}_t\bigl[X^{\bar{H}\bar{K}*}_t\bigr]\bigr)\\
			&\  +\bigl(\widetilde{C}^{i}-\widetilde{D}^{i} \bigl(\widetilde{\Sigma}^{\bar{H}\bar{K}}\bigr)^{-1}\bigl(\widetilde{L}^{\bar{H}\bar{K}}\bigr)^{\top}\bigr)\mathbb{E}^{\theta}_t\bigl[X^{\bar{H}\bar{K}*}_t\bigr]+\beta^{i}_{t}-\widetilde{D}^{i}\bigl(\widetilde{\Sigma}^{\bar{H}\bar{K}}\bigr)^{-1}\hbar^{\bar{H}\bar{K}}_{t}\bigr],\\
			k^{\bar{H}\bar{K},\ell j}_t=&\ \bigl(\Phi^{\bar{H}\bar{K}}_{(j,t)}-\Phi^{\bar{H}\bar{K}}_{(\ell,t)}\bigr)\bigl(X^{\bar{H}\bar{K}*}_t-\mathbb{E}^{\theta}_{t}\bigl[X^{\bar{H}\bar{K}*}_t\bigr]\bigr) +\bigl(\Pi^{\bar{H}\bar{K}}_{(j,t)}-\Pi^{\bar{H}\bar{K}}_{(\ell,t)}\bigr)\mathbb{E}^{\theta}_{t}\bigl[X^{\bar{H}\bar{K}*}_t\bigr]+\rho^{\bar{H}\bar{K},\ell j}_{t}.
		\end{aligned}\right.
	\end{equation*}By Theorem~\ref{the: main result definite},
	$(X^{\bar H\bar K*},u^{\bar H\bar K*},
	p^{\bar H\bar K},\zeta^{\bar H\bar K},k^{\bar H\bar K})$
	solves the transformed Hamiltonian system
	\eqref{equ: Ham for HK}.
	Applying the transformation \eqref{equ: Ham inve trans}
	with $H=\bar H$ and $K=\bar K$,
	we obtain a solution $(X^*,u^*,p,\zeta,k)$
	of \eqref{equ: Ham for com}
	by Lemma~\ref{lem: Ham of indefi}. The relations between $\Phi,\Pi$ and their transformed
	counterparts show that $(p,\zeta,k)$ satisfies \eqref{equ: p q k CMF}.
	
	By \eqref{equ: J HK and J}, the two cost functionals
	differ by a constant independent of the control.
	Therefore, $(X^*,u^*)$ is also the unique optimal pair
	of \textbf{Problem CMF}.
	
	Finally, \eqref{equ: J HK and J} gives
	\begin{equation*}
		\begin{aligned}
			J\bigl(\xi,\ell_0;u^*\bigr)
			&=
			J^{\bar{H}\bar{K}}
			\bigl(\xi,\ell_0;u^{\bar{H}\bar{K}*}\bigr)
			+\frac12\mathbb{E}\Bigl[
			|\xi-\mathbb{E}\xi|^2_{\bar{H}_{(\ell_0,0)}}
			+|\mathbb{E}\xi|^2_{\bar{K}_{(\ell_0,0)}}
			+\int_0^T\sum_{i=1}^{N}
			|\beta_t^i|^2_{\bar{H}_{(\theta_t,t)}}
			\,\mathrm{d}t
			\Bigr].
		\end{aligned}
	\end{equation*}
	Substituting the optimal cost of the transformed problem from
	Theorem~\ref{the: main result definite}
	and using \eqref{equ: Phi Pi HKbar}-\eqref{equ: L HKbar}, we obtain
	\begin{equation*}
		\begin{aligned}
			J\bigl(\xi,\ell_0; u^{*}\bigr)
			=&\ \frac{1}{2}\mathbb{E}\Bigl\{\int_{0}^{T}\Bigl[\sum_{i=1}^{N}\left|\beta^{i}_{t}\right|^{2}_{\Phi^{\bar{H}\bar{K}}_{(\theta_t,t)}}-\left|\hbar^{\bar{H}\bar{K}}_{t}\right|^{2}_{(\widetilde{\Sigma}^{\bar{H}\bar{K}}_{(\theta_t,t)})^{-1}}+2\langle \alpha_{t},\varphi^{\bar{H}\bar{K}}_{t}\rangle\Bigr]\mathrm{d}t\\
			&\  +\left|\xi-\mathbb{E}\xi\right|^{2}_{\Phi^{\bar{H}\bar{K}}_{(\ell_0,0)}}+\left|\mathbb{E}\xi\right|^{2}_{\Pi^{\bar{H}\bar{K}}_{(\ell_0,0)}}+2\langle \mathbb{E}\xi, \varphi^{\bar{H}\bar{K}}_{0}\rangle \Bigr\}\\
			&\ +\frac{1}{2}\mathbb{E}\Bigl\{\left|\xi-\mathbb{E}\xi\right|^{2}_{\bar{H}_{(\ell_0,0)}}+\left|\mathbb{E}\xi\right|^{2}_{\bar{K}_{(\ell_0,0)}}+\int_{0}^{T}\sum_{i=1}^{N}\left|\beta^{i}_{t}\right|^{2}_{\bar{H}_{(\theta_t,t)}}\mathrm{d}t\Bigr\}\\
			=&\  \frac{1}{2}\mathbb{E}\Bigl\{\int_{0}^{T}\Bigl[\sum_{i=1}^{N}\left|\beta^{i}_{t}\right|^{2}_{\Phi_{(\theta_t,t)}}-\left|\hbar_{t}\right|^{2}_{\widetilde{\Sigma}^{-1}_{(\theta_t,t)}}+2\langle \alpha_{t},\varphi_{t}\rangle \Bigr]\mathrm{d}t +\left|\xi-\mathbb{E}\xi\right|^{2}_{\Phi_{(\ell_0,0)}}+\left|\mathbb{E}\xi\right|^{2}_{\Pi_{(\ell_0,0)}}+2\langle \mathbb{E}\xi, \varphi_{0}\rangle \Bigr\}.
		\end{aligned}
	\end{equation*}
	This establishes the optimal cost formula \eqref{equ: exp cost of CMF}.
\end{proof}

\end{document}